\documentclass[10pt]{article}
\usepackage{amsmath, amsfonts, amsthm, amssymb,verbatim}
        \newtheorem{theorem}{Theorem}[section]
\newtheorem{definition}[theorem]{Definition}
        \newtheorem{proposition}[theorem]{Proposition}
        \newtheorem{lemma}[theorem]{Lemma}
        
        \newtheorem{remark}[theorem]{Remark}

\numberwithin{equation}{section}

    \newcommand{\IGNORE}[1]{}

    \newcommand{\D}{\mathcal{D}}
    \newcommand{\F}{\mathbf{F}}
    \newcommand{\N}{\mathbb{N}}
    \newcommand{\R}{\mathbb{R}}
    
    \newcommand{\W}{\mathrm{w}}

    \newcommand{\BA}{{\boldsymbol{a}}}
    \newcommand{\BD}{{\mathrm{b}}}
    
    \newcommand{\BZ}{{\boldsymbol{z}}}
    \newcommand{\CC}{{\mathcal{H}}}
    \newcommand{\CH}{{\bar{\CC}}}
    \newcommand{\CT}{{\bar{C}}}
    \newcommand{\DC}{\mathbf{D}}
    \newcommand{\DS}{\mathbf{d}}
    \newcommand{\GS}{\geqslant}
    \newcommand{\HA}{{\textstyle{\frac12}}}
    \newcommand{\KG}{{\textstyle\frac{\kappa}{\gamma-1}}}
    \newcommand{\LA}{\langle}
    \newcommand{\LS}{\leqslant}
    \newcommand{\OA}{\overline{a}}
    
    \newcommand{\OC}{\overline{c}}
    
    \newcommand{\OZ}{\overline{z}}
    \newcommand{\PV}{\mathrm{PV}}
    \newcommand{\RA}{\rangle}
    
    \newcommand{\SC}{{\mathbb{S}}}

    \newcommand{\UA}{\underline{a}}
    
    \newcommand{\UC}{\underline{c}}
    
    \newcommand{\UZ}{\underline{z}}

    \newcommand{\DST}{\displaystyle}
    \newcommand{\EPS}{\varepsilon}
    
    \newcommand{\IND}{\mathbf{1}}
    \newcommand{\LOC}{\mathrm{loc}}

    \newcommand{\OVA}{{\overline{A}}}
    
    \newcommand{\OVE}{{\overline{E}}}
    \newcommand{\OVM}{{\overline{M}}}
    \newcommand{\OVR}{{\overline{\rho}}}
    \newcommand{\OVU}{{\overline{u}}}
    \newcommand{\OVZ}{{\overline{\BZ}}}
    \newcommand{\TST}{\textstyle}
    \newcommand{\UNA}{{\underline{A}}}

    \newcommand{\MEAS}{M}
    \newcommand{\PROB}{\mathrm{Prob}}
    
    \newcommand{\WEAK}{
        \DOTSB\protect\relbar\protect\joinrel\rightharpoonup}

    \DeclareMathOperator*{\CI}{\mathrm{Ci}}
    
    \DeclareMathOperator*{\DIV}{\mathrm{div}}
    \DeclareMathOperator*{\SPT}{\mathrm{spt}}
    \DeclareMathOperator*{\CURL}{\mathrm{curl}}
    
    \DeclareMathOperator*{\ESUP}{\mathrm{ess\ sup}}
    \DeclareMathOperator*{\SIGN}{\mathrm{sign}}
 
\begin{document} 

\title
    {Finite energy solutions to the isentropic Euler equations with geometric effects}

\author
    {Philippe G. LeFloch$^1$ and Michael Westdickenberg$^2$} 
\renewcommand \today {} 
\date{\today} 

\maketitle

\footnotetext[1]{
     Laboratoire Jacques-Louis Lions \& Centre National de la Recherche Scientifique,
Universit\'e Pierre et Marie Curie (Paris 6), 4 Place Jussieu,  75252 Paris, France.
\\
E-mail : {pgLeFloch@gmail.com.} 
}

\footnotetext[2]{School of Mathematics,    Georgia Institute of Technology,
   686 Cherry Street,
     Atlanta, Georgia 30332-0160, USA. Email: mwest@math.gatech.edu. 
\textit{AMS Subject Classification.} 35L65, 76L05.
\textit{Keywords and phrases.} Isentropic Euler equations, spherical symmetry, global existence.
{\tt To cite this paper:} Jour. Math. Pures Appli. 88 (2007), 389--429. 
}

\begin{abstract}
Considering the isentropic Euler equations of compressible fluid
dynamics with geometric effects included, we establish the existence
of entropy solutions for a large class of initial data. We cover
fluid flows in a nozzle or in spherical symmetry when the origin
$r=0$ is included. These partial differential equations are
hyperbolic, but fail to be strictly hyperbolic when the fluid mass
density vanishes and vacuum is reached. Furthermore, when geometric
effects are taken into account, the $\sup$-norm of solutions can not
be controlled since there exist no invariant regions. To overcome
these difficulties and to establish an existence theory for
solutions with arbitrarily large amplitude, we search for solutions
with finite mass and total energy. Our strategy of proof takes
advantage of the particular structure of the Euler equations, and
leads to a versatile framework covering general compressible fluid
problems. We establish first higher-integrability estimates for the
mass density and the total energy. Next, we use arguments from the
theory of compensated compactness and Young measures, extended here
to sequences of solutions with finite mass and total energy. The
third ingredient of the proof is a characterization of the unbounded
support of entropy admissible Young measures. This requires the
study of singular products involving measures and principal values.
\end{abstract} 

\newpage 


\section{Introduction}\label{S:INTRO}

We are interested in the existence of entropy solutions to the Euler
equations for isentropic compressible fluids. Attention in the
literature has been so far restricted to {\em bounded} solutions and,
for this reason, current techniques apply to one-dimensional equations
or to simplified situations with symmetry only. Recall that the Euler
equations form a hyperbolic system of conservation laws; strict
hyperbolicity, however, fails when the fluid mass density
vanishes and vacuum is reached. This major difficulty for the analysis
was first dealt with by DiPerna \cite{DiP} using Tartar's method
of compensated compactness \cite{Tartar}.

When geometric effects are taken into account, the Euler equations
are no longer in a fully conservative form but consist of two
balance laws with variable coefficients. It is conceivable that due
to the interaction of characteristic waves and the geometry of the
problem, solutions may become unbounded at isolated points. For
spherically symmetric flows, for instance, the fluid can converge
towards the origin and waves can amplify nonlinearly, even if the
initial data was bounded pointwise. We are not aware of any result
showing that pointwise blow-up actually does occur. On the other
hand, there also seem to exist no method to establish boundedness in
full generality. In particular, the Conley-Chuey-Smoller principle
of invariant regions does not apply because the equations are not in
conservative form. Our objective is therefore to investigate the
isentropic Euler equations within a more general functional class:
We will only assume that solutions satisfy the natural bounds of
{\em finite mass and total energy.} The strategy we propose leads to
a versatile framework covering quite general compressible fluid
flows.

We are particularly interested in the case of spherically symmetric
flows where the origin $r=0$ is included in the domain, and of fluid
flows in a nozzle. Let us quickly recall the equations describing
these situations. We will assume that the nozzle is characterized by
a function $A=A(x)>0$ that determines its cross section at position
$x\in\R$. Then the isentropic Euler equations read
\begin{equation}
\begin{aligned}
    & \partial_t ( \rho A) + \partial_x ( \rho u A) = 0,
\\
    & \partial_t ( \rho u A) + \partial_x( \rho u^2 A)
        + A \partial_x P(\rho) = 0.
\end{aligned}
\label{E:CONS}
\end{equation}
The unknowns of this system are the density $\rho\GS 0$ and the
velocity $u$, which are functions of the independent variables
$(t,x) \in [0,\infty)\times\R$. The pressure $P(\rho)$ is related to
the internal energy $U(\rho)$ by the relation
$$
    P(\rho) = U'(\rho)\rho - U(\rho)
$$
for all $\rho\GS 0$. We restrict ourselves to polytropic perfect
gases, for which
$$
    U(\rho) = \KG\rho^\gamma
    \quad\text{and}\quad
    P(\rho)=\kappa\rho^\gamma.
$$
Here $\gamma>1$ is the adiabatic coefficient, and $\kappa :=
\theta^2/\gamma$ with $\theta := (\gamma-1)/2$ are constants. The
case of general pressure laws will be addressed in future work. The
first equation in \eqref{E:CONS} implies that the total mass is
conserved, thus
\begin{equation}
    M[\,\rho\,] := \int_\R \rho A \, dx
    \quad\text{is constant in time.}
\label{E:MD}
\end{equation}
The analogous statement for the momentum $\rho u A$ does not hold
because the momentum equation in general does not admit a conservative
form.

For spherically symmetric flows in $\R^d$, we have again equations
\eqref{E:CONS} with
$$
    A(x) := \omega_d x^{d-1}
    \quad\text{for all $x\in(0,\infty)$.}
$$
The constant $\omega_d>0$ denotes the volume of the unit sphere in
$\R^d$. Here the unknowns $(\rho,u)$ are defined for $(t,x)\in
[0,\infty) \times (0,\infty)$ and
$$
    M[\,\rho\,] := \int_{(0,\infty)} \rho A \,dx
    \quad\text{is constant in time.}
$$

In the following, we will cover both cases simultaneously by
considering the equations \eqref{E:CONS} with $A$ a continuously
differentiable function and
\begin{equation}
\begin{tabular}{|l|ll|}
    \hline\vphantom{\Big|}
    nozzle flow case
        & $\Omega:=\R$
        & \quad
            $A:\R \longrightarrow [\,\UNA, \OVA\,]$
\\
    \hline\vphantom{\Big|}
    spherical symmetry
        & $\Omega:= (0,\infty)$
        & \quad
            $A(x) := x^\alpha$
\\
    \hline
\end{tabular}
\label{E:EG.HYPO}
\end{equation}
Here, $\UNA < \OVA$ and $\alpha$ are positive constants. We also
require that
\begin{equation}
    (\partial_x A)_- \in L^1\cap L^\infty(\Omega),
\label{E:GEON}
\end{equation}
where $(b)_- := -\min\{b,0\}$ for all $b\in\R$. We refer the reader
to Sections~\ref{SS:HIGHER} and \ref{SS:EQUII} for further
explanation. Note that in the case of spherically symmetric flows
\eqref{E:GEON} is trivially satisfied since then $A$ is strictly
increasing. We also emphasize that for nozzle flows our arguments
can be adapted to work if assumption \eqref{E:GEON} is satisfied for
the positive part $(\partial_x A)_+$ instead. This is natural since
otherwise one direction would be favored, which would be unphysical.

It is easy to check that every {\em smooth} solution of
\eqref{E:CONS} admits an additional conservation law for the total
energy of the fluid
\begin{equation}
    \partial_t \Big( \big( \HA\rho u^2+U(\rho) \big) A\Big)
        + \partial_x \Big( \big( \HA \rho u^2 + Q(\rho) \big)
            u A \Big) = 0,
\label{E:CENER}
\end{equation}
where $Q(\rho) := U'(\rho)\rho$. The observation made earlier for
the mass equation applies again: the total energy associated with
{\em smooth} solutions of \eqref{E:CONS} is constant in time. For
{\em weak} solutions this equation
should not be imposed as an equality but as an inequality. In turn,
it is natural to require that for physically relevant weak
solutions
of \eqref{E:CONS}, the total energy
\begin{equation}
    E[\,\rho,u\,]
        := \int_\Omega \big( \HA\rho u^2+U(\rho) \big) A \,dx
    \quad\text{is nonincreasing in time.}
\label{E:ED}
\end{equation}

Our primary interest is about the Cauchy problem, so we impose the
condition
\begin{equation}
    \rho=\OVR,
    \qquad \rho u=\OVR \OVU
    \qquad \text{on $\{t=0\}\times\Omega$,}
\label{E:DATA}
\end{equation}
where $(\OVR, \OVU)$ is given initial data with finite mass and
total energy:
\begin{equation}
    M[\,\OVR\,] =: \OVM,
    \qquad E[\,\OVR,\OVU\,] =: \OVE,
    \qquad\text{with $\OVM, \OVE < \infty$.}
\label{E:FI}
\end{equation}

The selection of physically relevant solutions is based on a family
of entropy inequalities, which are defined as follows. For
$s \in\R$ and $(\rho,u)\in[0,\infty)\times\R$ introduce
the {\em entropy/entropy-flux kernels}
\begin{equation}
\begin{aligned}
&    \chi(s|\rho,u)
        :=  \Big( \rho^{2\theta}-(s-u)^2 \Big)_+^\lambda,
    \\
&    \sigma(s|\rho,u)
        :=  \Big( \theta s + (1-\theta) u \Big) \chi(s|\rho,u),
\end{aligned}
\label{E:CHISIG}
\end{equation}
where $\lambda := (3-\gamma)/2(\gamma-1)$ and $(b)_+ := \max\{b,0\}$
for all $b\in\R$. Observe that
$$
    \int_\R \begin{pmatrix}
            1 \\ s \\ \HA s^2
        \end{pmatrix} \Big( \chi(s|\rho,u), \sigma(s|\rho,u)
            \Big) \,ds
    = \begin{pmatrix}
        \rho & \rho u
\\
        \rho u & \rho u^2 + P(\rho)
\\
        \HA\rho u^2+U(\rho) &
            \big(\HA\rho u^2+Q(\rho) \big) u
    \end{pmatrix},
$$
which connects the Euler equations and the entropy/entropy-flux
kernels.

We will say that a function $\psi\in C^2(\R)$ is an {\em admissible
weight function} if it is convex and has subquadratic growth at
infinity. For all admissible weight functions $\psi$ we can
introduce the {\em entropy/entropy-flux pair}
\begin{equation}
    \Big( \eta_\psi(\rho,u), q_\psi(\rho,u) \Big)
        := \int_\R \psi(s) \,
            \Big( \chi(s|\rho,u), \sigma(s|\rho,u) \Big) \, ds,
\label{E:ETAPSI}
\end{equation}
and we impose the entropy inequalities
\begin{equation}
    \partial_t \Big( \eta_\psi(\rho,u) A \Big)
        + \partial_x \Big( q_\psi(\rho,u) A \Big)
        + \Big( \rho u \;\eta_{\psi,\rho}(\rho,u)
            - q_\psi(\rho,u) \Big) (\partial_x A)
    \LS 0
\label{E:ENTROPY}
\end{equation}
in the distribution sense. We use the notation $g_{,\rho} :=
\partial_\rho g$ for all functions $g$.

\begin{definition}
\label{D:SOLUTION} Let $(\OVR,\OVU)$ be given initial data with
finite mass and total energy. A pair of measurable functions
$(\rho,u): [0,\infty) \times \Omega \longrightarrow [0, \infty)
\times \R$ is called an entropy solution with finite mass and energy
(or a finite energy solution, for short) to the Cauchy problem
\eqref{E:CONS} \& \eqref{E:DATA} if the following is true:
\begin{enumerate}
\item The total mass is conserved in time: for almost every
(a.e.) $t$
$$
    M[\,\rho\,](t) = \OVM.
$$
\item The total energy is bounded in time: for a.e.\ $t$
$$
    E[\,\rho,u\,](t) \LS \OVE.
$$
\item The entropy inequalities \eqref{E:ENTROPY} are satisfied in
the distribution sense for all admissible weight functions $\psi$.
\item The initial data $(\OVR,\OVU)$ is
attained
in the distribution sense.
\end{enumerate}
\end{definition}

Clearly, the balance laws \eqref{E:CONS} follow from the entropy
inequality, by choosing $\psi$ to be constant or linear. Here is
our main result:

\begin{theorem}[Global Existence]
Consider the isentropic Euler equations \eqref{E:CONS} for a
polytropic perfect gas with adiabatic coefficient
$\gamma\in(1,5/3]$. Let the geometry be specified by
\eqref{E:EG.HYPO} \& \eqref{E:GEON}, where $\UNA < \OVA$ and
$\alpha$ are positive constants. Then, for any initial data $(\OVR,
\OVU)$ with finite mass and total energy, the Cauchy problem
\eqref{E:CONS} \& \eqref{E:DATA} admits a finite energy solution
$(\rho,u)$. \label{T:MAIN}
\end{theorem}

As we will show below, finite energy solutions have nonincreasing
total energy, so \eqref{E:ED} holds. But our estimates are not
strong enough to conclude that also a local energy balance is
satisfied (see Section~\ref{SS:EQUII} for further details). This is
the reason why only $\psi$ with subquadratic growth are considered
here. The local energy inequality can be recovered if we impose
higher-integrability for the initial data, as we will discuss in a
follow-up paper.

In the planar case, for which $A$ is constant, the existence of {\em
bounded} entropy solutions arising from bounded initial data was
first studied in pioneering work by DiPerna \cite{DiP}. His result
was generalized in \cite{Chen, CL, DCL1, DCL2, DiP, LPS, LPT}.
Existence of bounded solutions for the case of spherically symmetric
and nozzle flows were considered by Glimm and Chen \cite{CG}. To
avoid the difficulty of spherically symmetric solutions becoming
potentially unbounded, they constructed solutions outside a ball
around the origin only. A criterion for existence of bounded
solutions in the whole space (including the origin) was found by
Chen \cite{Chen2}: The inflow of the fluid towards the origin must
be below a certain threshold.

Our strategy to establish Theorem~\ref{T:MAIN} consists of two
parts. In Section~\ref{S:WEAKC} we first establish the existence of
measure-valued entropy solutions: We consider a sequence of bounded
approximate solutions $(\rho^n,u^n)$, obtained by suitably
truncating the unbounded initial data $(\OVR, \OVU)$ and then using
the existence results of \cite{CG}. We then prove the first key
observation that the approximate density $\rho^n$ enjoys
higher-integrability in space-time, i.e., we have
$$
    \text{$\rho^n \in L_\LOC^{\gamma+1}\big([0,\infty)\times
        \Omega \big)$ uniformly in $n$.}
$$
This fact is established by a commutator estimate, following a
strategy that was already used in \cite{DOW} in the context of
scalar conservation laws. A similar estimate was also derived in
\cite{Lions}. The second key observation made in
Section~\ref{S:WEAKC} is that also the total energy
$E[\,\rho^n,u^n\,]$ enjoys a higher integrability. The proof is
based on a bound for the entropy-flux, following the arguments in
\cite{LPT, Mi}. An alternative proof, which works for the planar
case only, is given in the Appendix. It relies on ``propagation of
equi-integrability'' for the total energy. The particular form of
the Euler equations and the freedom in choosing the weight function
$\psi$ in the definition of the entropy is essential here.

In Section~\ref{S:STRONG} we further analyze the structure of the
measure-valued solution. We show that the associated Young measure
$\nu_{(t,x)}$ is concentrated at a single point for almost every
$(t,x)$ and therefore conclude that the measure-valued solution is
actually a weak solution. This proves Theorem~\ref{T:MAIN}. To
achieve the Young measure reduction, we first apply compensated
compactness theory (see Tartar \cite{Tartar}) and derive the
well-known $\DIV$-$\CURL$-commutator relation. Then we determine the
support of the Young measure in the $(\rho,u)$-plane, for which we
must study singular products of distributions. Since we do not
require pointwise bounds on the solutions, we must also deal with
the difficulty that the support of the Young measure might be
unbounded.
\medskip

In the following, we denote by $C^k(B)$ the space of $k$-times
continuously differentiable functions, for suitable subsets
$B\subset R^N$. If $k=0$, then we simply write $C(B):=C^0(B)$. We
denote by $C_\BD(B)$ the space of bounded continuous functions,
whereas $C_0(B)$ is the closure of $\D(B)$ with respect to the
$\sup$-norm. Here, $\D(B)$ is the space of smooth functions with
compact support. The symbol $C^\alpha(B)$ with $\alpha\in(0,1)$ is
used for H\"{o}lder continuous functions.


\section{Weak convergence and measure-valued solutions}
\label{S:WEAKC}

In this section, we first construct a sequence of approximate
solutions $(\rho^n,u^n)$ to the isentropic Euler equations. These
functions are entropy solutions generated by compactly supported
bounded initial data. We then show the weak convergence of
approximate solutions to a measure-valued solution.

\subsection{Finite energy approximate solutions}\label{SS:APPROX}

In the spherically symmetric case, we need to remove the singularity
at the origin. We therefore introduce the modified geometry function
\begin{equation}
    A^n(x) := (x+1/n)^\alpha,
\label{E:EG.HYPOn}
\end{equation}
which converges uniformly to $A(x)=x^\alpha$ as
$n\rightarrow\infty$. The Cauchy problem associated to the function
$A^n$ is equivalent to a problem posed in the exterior of a ball of
radius $1/n$, for which existence of bounded entropy solution was
shown in \cite{CG}. In the case of nozzle flows we simply put
$A^n:=A$ for all $n$. Again we can use \cite{CG}. Let $M^n[\cdot]$
and $E^n[\cdot]$ denote the functionals defined in \eqref{E:MD} and
\eqref{E:ED}, with $A$ replaced by $A^n$. Given initial data $(\OVR,
\OVU)$ with $\OVR\GS 0$, we now consider a sequence of measurable
functions $(\OVR^n, \OVU^n)$ with $\OVR^n \GS 0$ that
\begin{enumerate}
\item are bounded and compactly supported in the closure
$\bar{\Omega}$;
\item converge in measure:
\begin{equation}
    \lim_{n\rightarrow\infty} (\OVR^n, \OVU^n) = (\OVR, \OVU);
\label{E:STRONG}
\end{equation}
\item have finite total mass $\OVM$:
\begin{equation}
    M^n[\,\OVR^n\,] = \OVM
    \quad\text{for all $n$;}
\label{E:INITMASS}
\end{equation}
\item have uniformly bounded total energy converging to $\OVE$:
\begin{equation}
    \sup_n E^n[\,\OVR^n,\OVU^n\,] \LS 2 \, \OVE,
    \qquad
    \lim_{n\rightarrow\infty} E^n[\,\OVR^n,\OVU^n\,] = \OVE.
\label{E:ENERCONV}
\end{equation}
\end{enumerate}
Clearly, it is possible to choose an approximating sequence
$(\OVR^n, \OVU^n)$ with the above properties, by first truncating
and mollifying the initial data $(\OVR, \OVU)$ and then
multiplying the density by a suitable constant to enforce
\eqref{E:INITMASS}.

Next, let $(\rho^n,u^n)$ be a sequence of entropy solutions of
\eqref{E:CONS} corresponding to the sequence of initial data
$(\OVR^n,\OVU^n)$. They have the following properties:
\begin{enumerate}
\item For any $n$ the entropy solution $(\rho^n,u^n)$ is bounded in
$L^\infty([0,\infty)\times\Omega)$ and has compact support in
space for all times $t\GS 0$.
\item The total mass is conserved in time: for a.e.\
$t$
\begin{equation}
    M^n[\,\rho^n\,](t) = M^n[\,\OVR^n\,].
\label{E:MASSCONS}
\end{equation}
\item The total energy is nonincreasing in time: for a.e.\
$t$
\begin{equation}
    E^n[\,\rho^n,u^n\,](t)
        \LS E^n[\,\OVR^n,\OVU^n\,].
\label{E:ENTRINEQ}
\end{equation}
\end{enumerate}
We will refer to a sequence of functions $(\rho^n, u^u)$ satisfying
the above conditions as a sequence of {\em finite energy
approximate solutions} of the Euler equations.

Our objective is to establish the strong pre-compactness of
$(\rho^n, u^n)$. To achieve this, we first derive a
higher-integrability property satisfied by the density $\rho^n$
uniformly in $n$. This will allow us to introduce a Young measure
representation for the limits of nonlinear functions of
$(\rho^n,u^n)$.

\subsection{Higher integrability of the mass density
variable}\label{SS:HIGHER}

We claim that for every $n$ there exists a function $h^n\colon
[0,\infty)\times\bar{\Omega} \longrightarrow \R$ that
\begin{enumerate}
\item has distributional derivatives
\begin{equation}
    \partial_t h^n = -\rho^n u^n A^n,
    \qquad \partial_x h^n = \rho^n A^n;
\label{E:PARTH}
\end{equation}
\item can be normalized so that
\begin{equation}
    0 \LS h^n \LS \OVM.
\label{E:NORM}
\end{equation}
In the spherically symmetric case, we may assume $h(t,0)=0$ for all
$t$.
\end{enumerate}
Note first that a function $h^n$ satisfying \eqref{E:PARTH} always
exists since the conservation law for $\rho$ precisely
says that the mixed second derivatives of $h^n$ commute. We see that
for almost every $t\GS 0$, the map $x\mapsto h^n(t,x)$ is absolutely
continuous and nondecreasing because the function $\rho^n A^n$ is
nonnegative.

Consider first the case of a nozzle, for which $\Omega=\R$. Since
the total mass is preserved we conclude that for a.e.\ $t\GS 0$ we
have the identity
\begin{equation}
    \lim_{x\rightarrow \infty} h^n(t,x)
        - \lim_{x\rightarrow -\infty} h^n(t,x) = \OVM.
\label{E:INTEGX}
\end{equation}
On the other hand, since for all fixed $t$ the functions
$(\rho^n,u^n)(t,\cdot)$ are compactly supported in $\R$ the first
identity in \eqref{E:PARTH} implies that
$$
    \lim_{x\rightarrow -\infty} h^n(t,x)
        = \lim_{x\rightarrow -\infty} h^n(0,x)
$$
for a.e.\ $t\GS 0$. Normalizing $h^n$ such that $\lim_{x
\rightarrow -\infty} h^n(0,x) = 0$, we get \eqref{E:NORM}.

Consider next the spherically symmetric case, for which  $\Omega =
(0,\infty)$. Then
\begin{equation}
    \lim_{x\rightarrow \infty} h^n(t,x)
        - \lim_{x\rightarrow 0} h^n(t,x) = \OVM
\label{E:INTEGX2}
\end{equation}
for a.e.\ $t\GS 0$. Since the momentum $\rho^n u^n A^n$ vanishes
at $x=0$, the first identity in \eqref{E:PARTH} implies that for
a.e.\ $t$ we obtain again
$$
    \lim_{x\rightarrow 0} h^n(t,x)
        = \lim_{x\rightarrow 0} h^n(0,x).
$$
Normalizing $h^n$ such that $\lim_{x \rightarrow 0} h^n(0,x) = 0$,
we again obtain \eqref{E:NORM}.

\begin{proposition}[Higher integrability]\label{P:HIGH}
Let $(\rho^n, u^n)$ be\ the finite energy approximate solutions
constructed in Subsection~\ref{SS:APPROX}, with geometry given by
\eqref{E:EG.HYPO} \& \eqref{E:GEON}. For any $T>0$ there exists a
constant $C>0$ such that
$$
    \sup_n \iint_{[0,T]\times\Omega} (\rho^n)^{\gamma+1}
        A^2 \,dx \,dt \LS C.
$$
\end{proposition}

{\em Proof.}
To simplify notation, we assume that in the spherically
symmetric case all functions are extended by zero for $x<0$. Recall
that we may assume the boundary condition $h^n(t,0)=0$ for all
$t$. Then \eqref{E:PARTH} holds in $[0,\infty)\times\R$.
\medskip

{\bf Step~1.} We will prove that $h^n$ is locally H\"{o}lder
continuous in both variables, with constants that are bounded
uniformly in $n$. The equi-continuity of $h^n$ in space
follows easily from \eqref{E:ENTRINEQ} and \eqref{E:PARTH}: Let $
K\subset\R$ be some compact subset. For all points $x_1,x_2\in K$
we
can then estimate
\begin{align*}
    & \ESUP_{t\GS 0} |h^n(t,x_2)-h^n(t,x_1)|
\\
    & \quad \LS \ESUP_{t\GS 0}  \int_{x_1}^{x_2} \rho^n A^n \,dx
\\
    & \quad \LS \ESUP_{t\GS 0} \bigg( \int_{x_1}^{x_2}
        (\rho^n)^\gamma A^n \,dx \bigg)^{\!1/\gamma}
    \bigg( \int_{x_1}^{x_2} A^n \,dx \bigg)^{\!(\gamma-1)/\gamma}.
\end{align*}
The first factor can be estimated by \eqref{E:ENTRINEQ} and
\eqref{E:ENERCONV}. We find
\begin{equation}
    \ESUP_{t\GS 0} |h^n(t,x_2)-h^n(t,x_1)|
        \LS C_1 |x_2 - x_1|^{(\gamma-1)/\gamma},
\label{E:SPACE}
\end{equation}
with $C_1>0$ some constant depending on $\OVE$ and
$\|A\|_{L^\infty(K)}$ (recall \eqref{E:EG.HYPOn}).

To prove the equi-continuity in time we first fix a mollifier
$\varphi_\delta$ with the standard properties $\varphi_\delta \geq
0$, $\int \varphi_\delta \, dx =1$, and $\SPT \varphi_\delta
\subset (-\delta, \delta)$. The parameter $\delta>0$ will be chosen
later on. We then deduce from \eqref{E:SPACE} that for all $x\in K$
\begin{align*}
    & \vphantom{\int_\R}
        \ESUP_{t\GS 0} \bigg| \bigg( \int_\R
            \varphi_\delta(x-y) h^n(t,y) \,dy
                \bigg)-h^n(t,x) \bigg|
\\
    & \quad \LS C_1 \int_\R \varphi_\delta(x-y)
        |x-y|^{(\gamma-1)/\gamma} \,dy
\\
    & \vphantom{\int_\R}
        \quad \LS C_1 \delta^{(\gamma-1)/\gamma}.
\end{align*}
For any $t_1,t_2\GS 0$ and $x\in\R$ we therefore obtain
\begin{align}
    & |h^n(t_2,x)-h^n(t_1,x)|
\nonumber\\
    & \quad \LS 2C_1\delta^{(\gamma-1)/\gamma}
            + \bigg| \int_\R \varphi_\delta(x-y)\,
                \big( h^n(t_2,y)-h^n(t_1,y) \big) \,dy \bigg|
\nonumber\\
    & \quad = 2C_1\delta^{(\gamma-1)/\gamma}
            + \bigg| \int_{t_1}^{t_2} \int_\R \varphi_\delta(x-y)\,
                (\rho^n u^n)(t,y) A^n(y) \,dy \,dt \bigg|.
\label{E:TIMEC}
\end{align}
Now note that the energy bound \eqref{E:ENTRINEQ} implies the
estimate
\begin{align}
    & \ESUP_{t\GS 0} \int_\R
            |\rho^n u^n|^{2\gamma/(\gamma+1)} A^n \,dx
\nonumber\\
    & \quad \LS \ESUP_{t\GS 0} \bigg( \int_\R (\rho^n)^\gamma
            A^n \,dx \bigg)^{\!1/(\gamma+1)}
        \bigg( \int_\R \rho^n (u^n)^2
            A^n \,dx \bigg)^{\!\gamma/(\gamma+1)}
\nonumber\\
    & \vphantom{\int_\R}
        \quad \LS C_2,
\label{E:MOM}
\end{align}
with $C_2>0$ some constant depending on \eqref{E:ENERCONV}. Using
this in \eqref{E:TIMEC} and optimizing in $\delta$, we arrive at the
following estimate: for any $t_1,t_2\GS 0$
\begin{align*}
    & \ESUP_{x\in\R} |h^n(t_2,x)-h^n(t_1,x)|
\\
    & \quad \LS 2C_1\delta^{(\gamma-1)/\gamma}
        + C_2^{(\gamma+1)/2\gamma}\|\varphi\|_{L^\infty(\R)}
            \delta^{-(\gamma+1)/2\gamma}|t_1-t_2|
\\
    & \vphantom{\int_\R}
        \quad \LS C_3 |t_1-t_2|^{2(\gamma-1)/(3\gamma-1)}
\end{align*}
for some constant $C_3>0$. This establishes the first part of the
proposition.
\medskip

{\bf Step~2.} Let $\varphi_\EPS$ be a standard mollifier in
$\R^2$ and, after extending $h^n$ by zero to all of $\R^2$, define
the smooth function $h^n_\EPS := h^n\star\varphi_\EPS$. Then the
following identity is true in the distribution sense in
$[0,\infty)\times\R$:
\begin{align}
    & \partial_t \Big(  \rho^n u^n A^n \; h^n_\EPS \Big)
       + \partial_x \Big( \rho^n (u^n)^2 A^n \; h^n_\EPS \Big)
       + A^n \partial_x \Big( P(\rho^n) \; h^n_\EPS \Big)
\\
    & \quad
        = \bigg\{ \partial_t (\rho^n u^n A^n)
        + \partial_x \big( \rho^n (u^n)^2 A^n \big)
        + A^n \partial_x P(\rho^n) \bigg\}
            \; h^n_\EPS
\\
    & \qquad
        + \bigg\{ \rho^n u^n A^n \; (\partial_t h^n_\EPS)
        + \big( \rho^n (u^n)^2 + P(\rho^n) \big) A^n
            \; (\partial_x h^n_\EPS) \bigg\}.
\end{align}
The first term on the right-hand side vanishes in view of the
momentum conservation law satisfied by $(\rho^n,u^n)$. As
$\EPS\rightarrow 0$, we have $h^n_\EPS\to h^n$ uniformly on compact
sets because $h^n$ is equi-continuous by Proposition~\ref{P:HIGH}.

On the other hand, we have $\partial_t h^n_\EPS \to \partial_t h^n$
and $\partial_x h^n_\EPS \to \partial_x h^n$ in $L^1_\LOC([0,\infty)
\times\R)$. By boundedness of $(\rho^n,u^n)$ and \eqref{E:PARTH}, we
find that in distributional sense
\begin{align}
    P(\rho^n) \rho^n (A^n)^2
        & = \partial_t \Big( \rho^n u^n A^n \; h^n \Big)
        + \partial_x \Big( \big( \rho^n (u^n)^2 + P(\rho^n) \big)
            A^n \; h^n \Big)
\nonumber\\
        & \quad
            - h^n P(\rho^n) \; (\partial_x A^n).
\label{E:MORE}
\end{align}
We test \eqref{E:MORE} against a monotone sequence of functions
$\zeta_k\in \D([0,\infty)\times\bar{\Omega})$ with $0\LS\zeta_k\LS
1$ and $\zeta_k\to \mathbf{1}_{[0,T]\times\Omega}$ for some $T>0$.
Note that \eqref{E:ENTRINEQ} implies
\begin{align*}
    & \ESUP_{t\GS 0} \int_\R |\rho^n u^n| A^n \,dx
\\
    & \quad
        \LS \ESUP_{t\GS 0} \bigg( \int_\R \rho^n
                A^n \,dx \bigg)^{\!1/2}
            \bigg( \int_\R \rho^n (u^n)^2
                A^n \,dx \bigg)^{\!1/2},
\end{align*}
which can be estimated against $\sqrt{2\OVM\OVE}$. Since $(\rho^n,
u^n)$ has compact support in $x$ and since $h^n\GS 0$ is uniformly
bounded by $\OVM$, we obtain that for all $n$
\begin{equation}
    \iint_{[0,T]\times\Omega} (\rho^n)^{\gamma+1}
            (A^n)^2 \,dx \,dt
        \LS 2\OVM\, {\TST\sqrt{2\OVM\OVE}}
            + T\, \OVM\OVE\,
                \|(\partial_x A)_-\|_{L^\infty(\R)}.
\label{E:UNICON}
\end{equation}
For the spherically symmetric case we used the fact that $h^n$
vanishes at the origin, so the $x$-derivative on the left-hand side
of \eqref{E:MORE} does not contribute. We have $\partial_x A^n
\longrightarrow \partial_x A$ because of \eqref{E:EG.HYPOn} and
$(\partial_x A)_- = 0$. Therefore the second term in the estimate
\eqref{E:UNICON} vanishes in that case. Finally, note that $A\LS
A^n$ for all $n$, which proves the proposition in the case of
spherical symmetry. For nozzle flows we defined $A^n:=A$ for all
$n$, so there is nothing more to prove. Note that by normalizing the
function $h^n$ such that $-\OVM \LS h^n \LS 0$, we can also obtain
\eqref{E:UNICON} with $(\partial_x A)_-$ replaced by the positive
part of the gradient. \qed

Note that for any compact subset $K\subset[0,\infty)\times\Omega$
the function $A^2$ can be estimated uniformly from above and below.
In view of \eqref{E:EG.HYPO} this is obvious for the nozzle flow
case. For the case of spherically symmetric flows, observe that the
compact set $K$ is bounded away from the origin because $\Omega =
(0,\infty)$ is an open set. Proposition~\ref{P:HIGH} therefore
implies that
$$
    \text{$\rho^n \in L_\LOC^{\gamma+1}\big([0,\infty)\times
        \Omega \big)$ uniformly in $n$.}
$$
%


\subsection{Young measures based on energy bounds}\label{SS:YOUNG}

It will be convenient to work with the Riemann invariants
$(\OZ,\UZ)$ associated with \eqref{E:CONS}, rather than with the
physical variables $(\rho, u)$. For simplicity of notation, we will
consistently denote pairs of numbers such as $(\OZ,\UZ)$ by the
corresponding bold symbol $\BZ:=(\OZ,\UZ)$. We have
\begin{equation}
    \OZ(\rho,u) = u+\rho^\theta,
    \qquad
    \UZ(\rho,u) = u-\rho^\theta,
\label{E:RIEMANN}
\end{equation}
which is equivalent to
\begin{equation}
    \rho(\BZ) = \bigg(\frac{\OZ-\UZ}2\bigg)^{\!1/\theta},
    \qquad
    u(\BZ) = \frac{\OZ+\UZ}2.
\label{RIEMANN2}
\end{equation}
We consider entropies/entropy-fluxes as functions of $(\rho,u)$ or
$\BZ$, respectively.

We now define $H := \{ \BA\in\R^2\colon \OA>\UA\}$, and we will
tacitly assume that all functions in $\D(H)$ or $C_0(H)$ are
extended by zero to the closure $\bar{H}$, if necessary. Consider
then the following space of bounded continuous functions
\begin{align*}
    \CT(H) := \bigg\{ \; \varphi\in C(\bar{H})\colon
        & \text{the function $\varphi$ is constant in
            $\{\BA\in\R^2\colon \OA=\UA\}$ and}
\\[-2ex]
        & \text{the map $\DST\Big(\BA \mapsto
            \lim_{s\rightarrow\infty} \varphi(s\BA)\Big)$
                belongs to $C(S^1\cap\bar{H})$} \;\; \bigg\},
\end{align*}
where $S^1\subset\R^2$ denotes the sphere. This space allows us to
deal with the two difficulties of the problem under consideration:
at the vacuum and in the large. Observe that $\CT(H)$ has a ring
structure and is complete with respect to the $\sup$-norm.
Therefore, there exists a compactification $\CH$ of $H$ such that
$\CT(H)$ is isomorphic to the space $C(\CH)$. We refer the reader to
\cite{Roubicek, Rudin}. For simplicity, we will not distinguish
between functions in $\CT(H)$ and in $C(\CH)$.

The topology of $\CH$ is the weak-$\star$ topology induced
by $C(\CH)$: the sequence of points $\BA_n\in\CH$ converges to
$\BA \in\CH$ as $n\rightarrow\infty$ if and only if
$$
    \lim_{n\rightarrow\infty} \varphi(\BA_n) = \varphi(\BA)
    \quad\text{ for all $\varphi\in C(\CH)$.}
$$
In $H\subset\CH$ this weak-$\star$ topology is consistent with the
Euclidean topology, and thus $\CH$ is separable.
Moreover, the space $\CH$ is metrizable since $\CT(H)$ is separable
and separates points in $H$ (see Proposition~1.5.3 of
\cite{Roubicek} and Section~3.8 of \cite{Rudin}). On the other
hand, we emphasize the fact that the topology above does not
distinguish points in the compactification of the diagonal
$\{\BA\in\R^2\colon \OA=\UA\}$. In that sense, all points in the
vacuum are equivalent. We denote by $V$ the compactification of
$\{\BA\in\R^2\colon \OA=\UA\}$, and we define $\CC := H\cup V$.

We need the following result (see Theorem~2.4 of \cite{AM}).

\begin{theorem}[Young measures]\label{T:YOUNG}
Given any sequence of measurable functions $\BZ^n\colon
[0,\infty)\times \Omega \rightarrow \CH$ there exists a subsequence
(still labeled $\BZ^n$) and a function $\nu\in
L^\infty_\W([0,\infty) \times \Omega,\PROB(\CH))$ (that is, a
weakly-$\star$ measurable map from $[0,\infty)\times \Omega$ into
the space of probability measures on $\CH$), such that
$$
    \varphi(\BZ^n) \WEAK
        \int_\CH \varphi(\BA) \,\nu(d\BA)
    \quad\text{ weakly-$\star$ in
         $L^\infty\big([0,\infty)\times\Omega\big)$
         for all $\varphi\in C(\CH)$.}
$$
The functions $\BZ^n$ converge in measure to $\BZ\colon
[0,\infty)\times \Omega \rightarrow \CH$ if and only if
$$
    \nu_{(t,x)} = \delta_{\BZ(t,x)}
    \quad\text{for a.e.\ $(t,x)$.}
$$
\end{theorem}

We will use Young measures to represent limits of certain nonlinear
functions of $(\BZ^n)$ that may be unbounded. Let us introduce the
weight function
$$
    W(\BA)
        := 1 + \rho(\BA)^{\gamma+1}
    \quad\text{for all $\BA\in H$.}
$$
%

\begin{proposition}\label{P:FIRST} Consider the sequence of Riemann
invariants $(\BZ^n)$ associated with the sequence of finite energy
approximate solutions $(\rho^n,u^n)$ of Subsection~\ref{SS:APPROX}.
Let $\nu$ be a Young measure generated by (a subsequence of)
$(\BZ^n)$. Then for almost every $(t,x)\in[0,\infty)\times\Omega$ we
have that
\begin{equation}
    \nu_{(t,x)}\in\PROB(\CC),
    \qquad
    \int_\CC W(\BA) \,\nu_{(t,x)}(d\BA) < \infty.
\label{E:NUPROP}
\end{equation}
For any $\varphi = \varphi_0W$ with $\varphi_0 \in C_0(H)$ it holds
\begin{equation}
    \varphi(\BZ^n) \WEAK \LA\varphi\RA
         := \int_\CC \varphi(\BA) \,\nu(d\BA)
    \quad\text{weakly in
        $L^1_\LOC\big([0,\infty)\times\Omega\big)$.}
\label{E:RES}
\end{equation}
\end{proposition}

\begin{remark}\label{R:CONME}
The first statement in \eqref{E:NUPROP} means that $\nu_{(t,x)}$ is
supported in $H \cup V$ only instead of $\CH$. Note that in
\eqref{E:RES} we consider local convergence in the {\em open} set
$\Omega$. For the spherically symmetric case, this means convergence
away from the origin. A slightly more precise statement is
$$
    \varphi(\BZ^n) (A^n)^2 \WEAK \LA\varphi\RA A^2
    \quad\text{weakly in
        $L^1_\LOC\big([0,\infty)\times\bar{\Omega}\big)$}
$$
for all $\varphi = \varphi_0W$ with $\varphi_0 \in C_0(H)$.
Recall that $A^n$ converges uniformly to $A$.
\end{remark}


{\em Proof.}  We proceed in three steps.
\medskip

{\bf Step~1.} Let $\bar{B}_r(0)$ be the closed ball with radius $r$.
Fix a radial test function $\varphi\in C(\bar{H})$ with
$0\LS\varphi\LS 1$, such that $\varphi = 1$ in $\bar{H}\cap
\bar{B}_1(0)$ and $\varphi=0$ for $\bar{H} \setminus B_2(0)$. Let
$\varphi_R  :=   \varphi(\cdot/R)$ and $\Phi_R  :=  1-\varphi_R$ for
all $R>0$. Choose $\phi\in C(S^1\cap\bar{H})$ with $0\LS\phi\LS 1$
and compactly supported in $S^1 \cap H$, and extend $\phi$ as a
homogeneous function of degree zero to $\bar{H}\setminus\{0\}$. Then
$\phi\Phi_R \in \CT(H)$, so it can be identified with a function in
$C(\CH)$. Now Theorem~\ref{T:YOUNG} applies, and we obtain that for
any compact set $K\subset [0,\infty)\times\Omega$
\begin{align*}
    \iint_K \bigg( \int_\CH \phi(\BA) \Phi_R(\BA)
            \,\nu_{(t,x)}(d\BA) \bigg) \,dx\,dt
        & = \lim_{n\rightarrow\infty} \iint_K \phi(\BZ^n)
            \Phi_R(\BZ^n) \,dx\,dt
\\
        & \LS \sup_n \big| \big\{ \OZ^n-\UZ^n \GS c_\phi R
            \big\}\cap K \big|,
\end{align*}
where the constant $c_\phi>0$ depends on the support of $\phi$.
Hence, we get
\begin{align*}
    & \iint_K \bigg( \int_\CH \phi(\BA) \Phi_R(\BA)
            \,\nu_{(t,x)}(d\BA) \bigg) \,dx\,dt
\\
    & \quad
        \LS \frac{1}{1+\big(\frac{c_\phi R}{2}
            \big)^{(\gamma+1)/\theta}}\;
            \sup_n \iint_K W(\BZ^n) \,dx\,dt
        \longrightarrow 0
    \quad\text{ as $R\rightarrow\infty$.}
\end{align*}
Note that $W(\BZ^n)$ is uniformly bounded in $L^1(K)$ because of
Proposition~\ref{P:HIGH} and our assumptions on $A^n$ and $K$.
Since $\phi$ and $K$ were arbitrary, we conclude that $\nu$ is
supported in $H$ and the vacuum, thus $\nu_{(t,x)} \in \PROB(\CC)$
a.e.
\medskip

{\bf Step~2.} Consider a monotone sequence of $\phi_k\in\D(H)$ with
$0\LS\phi_k\LS 1$ and $\phi_k\to 1$ pointwise as $k\to\infty$. For
any $K\subset[0,\infty)\times \R$ compact we have
$$
    \iint_K \LA W\RA \, dx \,dt
        = \lim_{k\rightarrow\infty}
            \iint_K \LA \phi_k W\RA \, dx \,dt,
$$
by monotone convergence. On the other hand, Theorem~\ref{T:YOUNG}
yields
\begin{align*}
    \iint_K \LA \phi_k W\RA \, dx \,dt
    & = \lim_{n\rightarrow\infty} \iint_K
            \phi_k(\BZ^n) W(\BZ^n) \, dx\,dt
    \\
    & \LS \sup_n \iint_K W(\BZ^n) \, dx\,dt,
\end{align*}
which is finite by Proposition~\ref{P:HIGH} and by choice of $A^n$
and $K$.
\medskip

{\bf Step~3.} Let now $\varphi_0\in C_0(H)$ and choose a sequence
of functions $\varphi_k\in \D(H)$ with $\varphi_k \rightarrow
\varphi_0$ in the $\sup$-norm as $k\rightarrow \infty$. For any
$K\subset[0,\infty)\times \Omega$ compact and $\zeta \in
C_\BD([0,\infty)\times\Omega)$ and by setting $\varphi =
\varphi_0 W$, we can then estimate
\begin{align*}
    & \bigg| \iint_K \LA\varphi\RA\zeta \,dx\,dt
        - \iint_K \varphi(\BZ^n)\zeta \,dx\,dt \bigg|
\\
    & \quad
        \LS \|\varphi_k-\varphi_0\|_{L^\infty(H)}
            \|\zeta\|_{L^\infty(K)} \bigg(
                \iint_K \LA W\RA \,dx \,dt
                + \sup_n \iint_K W(\BZ^n) \, dx\,dt \bigg)
\\
    & \qquad
        + \bigg| \iint_K \LA\varphi_k W\RA \zeta \, dx\,dt
            -\iint_K \varphi_k(\BZ^n) W(\BZ^n) \zeta \, dx\,dt
                \bigg|
        \longrightarrow 0
    \quad\text{as $k,n\rightarrow\infty$.}
\end{align*}
Indeed, the first term on the right-hand side vanishes as
$k\rightarrow\infty$, by choice of $\varphi_k$ and in view of Step~2
and Proposition~\ref{P:HIGH}. The second term vanishes for any fixed
$k$ as $n\rightarrow\infty$, by Theorem~\ref{T:YOUNG}. This
completes the proof. \qed


\subsection{Measure-valued solutions}\label{SS:MEASURE}

Recall first  that in the seminal work \cite{LPT} the authors
introduced the kinetic formulation for the isentropic Euler
equations. They showed that for bounded entropy solutions, the
requirement that the inequality \eqref{E:ENTROPY} holds for a
sufficiently large class of admissible weight functions $\psi$, can
be reformulated in terms of a single kinetic equation with suitable
source term. This result can be generalized to the isentropic Euler
equations with geometric effect as follows: Let $(\chi,\sigma)$ be
the entropy/entropy-flux kernels introduced in \eqref{E:CHISIG}.
Then the pair of functions $(\rho,u)$ is a finite energy solution
of
\eqref{E:CONS} \& \eqref{E:DATA} if and only if there exists a
nonnegative bounded measure $\mu$ depending on
$(t,x)\in[0,\infty)\times\Omega$ and $s\in\R$ such that in the
distribution sense in $([0,\infty)\times\Omega)\times\R$ we have
\begin{align}
    & \partial_t \Big( \chi(\cdot|\rho,u) A \Big)
        + \partial_x \Big( \sigma(\cdot|\rho,u) A \Big)
        + \Big( \rho u\; \chi_{,\rho}(\cdot|\rho,u)
            -\sigma(\cdot|\rho,u) \Big) (\partial_x A)
\nonumber\\
    & \quad
        = -\partial^2_s(A\mu).
\label{E:KINETIC}
\end{align}
Recall that a finite energy solution satisfies the entropy
inequality \eqref{E:ENTROPY} for a large class of convex weights
$\psi$. The proof of this kinetic formulation follows closely the
one given in \cite{LPT} for the planar case (see also \cite{Mi} for
spherically symmetric flows), and we refer the reader to the
literature for further details. The measure $\mu$ captures the
entropy dissipation. It can be bounded as
\begin{equation}
    \iint_{[0,\infty)\times\Omega} \int_\R A(x)\, \mu(ds,dx,dt)
        \LS \int_\R \Big( \HA\OVR \OVU^2
            + U(\OVR) \Big) A \,dx.
\label{E:BDMEAS}
\end{equation}
A similar kinetic formulation can be derived for the sequence of
finite energy approximate solutions $(\rho^n,u^n)$ constructed in
Section~\ref{SS:APPROX}.

We are going to show now that a suitable subsequence of
$(\rho^n,u^n)$ converges to a measure-valued solution of the
isentropic Euler equations. In slight abuse of notation, we will
occasionally consider the entropy/entropy-flux kernels $(\chi,
\sigma)$ as functions of the Riemann invariants $\BZ$ instead of
$(\rho,u)$: We write
\begin{align*}
    & \chi(s|\BZ)
        :=  \Big( (\OZ-s)(s-\UZ) \Big)^\lambda_+,
\\
    & \sigma(s|\BZ)
         :=  \bigg( \theta s+(1-\theta) \frac{\OZ+\UZ}{2} \bigg)\,
            \chi(s|\BZ)
\end{align*}
for $s\in\R$, which is consistent with \eqref{E:CHISIG} (see
\eqref{E:RIEMANN}).

We need the following two observations.

\begin{lemma}
\label{L:PSIIND}
Assume that the sequence $(\rho^n,u^n)$ of finite
energy approximations constructed in Section~\ref{SS:APPROX}
generates a Young measure $\nu$ as explained in
Proposition~\ref{P:FIRST}. Let $(\BZ^n)$ be the Riemann
invariants associated with $(\rho^n,u^n)$. For any $\psi\in\D(\R)$,
the pair $(\eta_\psi, q_\psi)$ defined by \eqref{E:ETAPSI} then
satisfies
\begin{equation}
    \begin{aligned}
        \eta_\psi(\BZ^n) &\WEAK \LA\eta_\psi\RA
\\
        q_\psi(\BZ^n) &\WEAK \LA q_\psi\RA
    \end{aligned}
    \quad\text{weakly in
        $L^{\gamma+1}_\LOC\big([0,\infty)\times\Omega\big)$.}
\label{E:WEAKC}
\end{equation}
We also have
\begin{equation}
    (\rho u\; \eta_{\psi,\rho})(\BZ^n)
        \WEAK \LA \rho u\; \eta_{\psi,\rho} \RA
    \quad\text{weakly in
        $L^2_\LOC\big([0,\infty)\times\Omega\big)$.}
\label{E:WEAKC2}
\end{equation}
Moreover, if $\eta_{\psi'}$ is defined as in \eqref{E:ETAPSI} for
some $\psi'\in\D(\R)$, then
$$
    \begin{aligned}
        \eta_\psi(\BZ^n) \eta_{\psi'}(\BZ^n)
            &\WEAK \LA\eta_\psi \eta_{\psi'}\RA
\\
        q_\psi(\BZ^n) \eta_{\psi'}(\BZ^n)
            &\WEAK \LA q_\psi \eta_{\psi'}\RA
    \end{aligned}
    \quad\text{weakly in
        $L^1_\LOC\big([0,\infty)\times\Omega\big)$.}
$$
\end{lemma}

{\em Proof.}
A straightforward change of variables shows that $\eta_\psi$ is
given by
\begin{equation}
    \eta_\psi(\BA)
        = \rho(\BA) \int_{-1}^1 \psi\Big( u(\BA)+t\rho(\BA)^\theta
            \Big) \, (1-t^2)^\lambda \,dt,
\label{E:ERT}
\end{equation}
so clearly $\BA\mapsto\eta_\psi(\BA)$ is a continuous function.
Suppose that the support $\SPT\psi$ of the function $\psi$ is
included in an interval $[\UC,\OC]$. Then we have
\begin{equation}
    |\eta_\psi(\BA)|
        \LS C \, \IND_{\{\UC\LS\UA\}}
            \IND_{\{\OA\LS\OC\}}
        \begin{cases}
            \;\rho(\BA)\quad
                & \text{for $\OA-\UA$ small,}
\\
            \;\rho(\BA)^{2\lambda\theta}
                & \text{for $\OA-\UA$ large,}
        \end{cases}
\label{E:ETAEST}
\end{equation}
with $C>0$ a constant depending on $\psi$ and $\lambda$. Indeed,
note that $\lambda > 0$ for $\gamma\in(1,3)$, which implies that the
map $t\mapsto (1-t^2)^\lambda$ is integrable on $[-1,1]$. The
behavior for small $\OA-\UA$ then follows immediately. For large
$\OA-\UA$, the $s$-integral in \eqref{E:ERT} is restricted to an
interval of length $(\OC-\UC) / \rho(\BA)^\theta$. This implies that
the integral in \eqref{E:ERT} is bounded above by a constant times
$1/\rho(\BA)^{\theta}$. Since $1-\theta=2\lambda\theta$, the
asymptotic behavior in \eqref{E:ETAEST} follows. We conclude that
$$
    \eta_\psi W^{-1} \in C_0(H)
    \quad\text{and}\quad
    \eta_\psi \eta_{\psi'} W^{-1} \in C_0(H)
$$
(since $4\lambda\theta<\gamma+1$ if $\gamma>1$), and by
Proposition~\ref{P:FIRST}
\begin{equation}
    \begin{aligned}
        \eta_\psi(\BZ^n)
            &\WEAK \LA\eta_\psi\RA
\\
        \eta_\psi(\BZ^n) \eta_{\psi'}(\BZ^n)
            &\WEAK \LA\eta_\psi \eta_{\psi'}\RA
    \end{aligned}
    \quad\text{ weakly in
        $L^1_\LOC\big([0,\infty)\times\Omega\big)$.}
\label{E:WELK}
\end{equation}
We also have $|\eta_\psi(\BA)|^{\gamma+1} \LS C \,
W(\BA)$ for all $\BA\in H$ and some constant $C>0$. Therefore
\eqref{E:WELK} can be improved to \eqref{E:WEAKC}, in view of
Proposition~\ref{P:HIGH}.

For $q_\psi$ we can argue in a similar way, using the bound
\begin{align}
    |q_\psi(\BA)|
        & \LS \max\big( |\OA|,|\UA|\big) \, |\eta_\psi(\BA)|
\nonumber\\
        & \LS \Big( \max\big\{|\OC|,|\UC|\big\} + (\OA-\UA) \Big)
            \, |\eta_\psi(\BA)|
    \quad\text{for all $\BA\in H$.}
\label{E:QEST}
\end{align}
We have $q_\psi W^{-1} \in C_0(H)$ and $q_\psi \eta_{\psi'} W^{-1}
\in C_0(H)$ (since $(4\lambda+1)\theta<\gamma+1$), and
$|q_\psi(\BA)|^{\gamma+1} \LS CW(\BA)$ for all $\BA\in H$ and some
constant $C>0$.

The statement in \eqref{E:WEAKC2} follows analogously. We use
the identity
\begin{align*}
    (\rho u\; \eta_{\psi,\rho})(\BA)
        & = u(\BA) \int_\R \psi(s) \chi(s|\BA) \,ds
\nonumber\\
        & \quad
            + \theta u(\BA) \int_\R \psi'(s)\big(s-u(\BA)\big)
                \chi(s|\BA) \,ds,
\end{align*}
and then proceed as in \eqref{E:QEST}. Note that
$2(\lambda+1)\theta = (\gamma+1)/2$.
\qed

We now establish strong convergence of the approximate initial
data.

\begin{lemma}
\label{L:CONID}
For any smooth weight function $\psi$ with at most
quadratic growth at infinity, let the entropy $\eta_\psi$ be
defined
by \eqref{E:ETAPSI}. Then we have
$$
    \eta_\psi(\OVR^n,\OVU^n)
        \longrightarrow \eta_\psi(\OVR,
            \OVU)
    \quad\text{strongly in $L^1_\LOC(\Omega)$.}
$$
\end{lemma}

{\em Proof.}  By assumption~\eqref{E:STRONG}, we have $(\OVR^n,
\OVU^n) \longrightarrow (\OVR,\OVU)$ in measure. It therefore
suffices to show equi-integrability of $\eta_\psi (\OVR^n, \OVU^n)$
locally. We choose a function $\varphi\in\D(\R)$ with
$0\LS\varphi\LS 1$, such that $\varphi(s)=1$ for $|s|\LS 1$ and
$\varphi(s)=0$ for $|s|\GS 2$. Define $\varphi_R :=
\varphi(\cdot/R)$ and $\Phi_R  := 1-\varphi_R$, and fix some
$K\subset\Omega$ compact. We will show that for all $\EPS>0$ there
exist numbers $N,R>0$ with
\begin{equation}
    \sup_{n\GS N} \iint_{K\times\R} s^2\Phi_R(s)\,
        \chi(s|\OVZ^n) \,ds \,dx \LS \EPS.
\label{E:CLA}
\end{equation}
Indeed, we can decompose
\begin{align}
    & \iint_{K\times\R} s^2\Phi_R(s)\, \chi(s|\OVZ^n) \,ds \,dx
\nonumber\\
    & \quad
        = \bigg( \iint_{K\times\R} s^2\chi(s|\OVZ^n) \,ds \,dx
            -\iint_{K\times\R} s^2\chi(s|\OVZ) \,ds \,dx \bigg)
\nonumber\\
    & \qquad
        - \bigg( \iint_{K\times\R} s^2\varphi_R(s)\,
                \chi(s|\OVZ^n) \,ds \,dx
            - \iint_{K\times\R} s^2\varphi_R(s) \,
                \chi(s|\OVZ) \,ds \,dx \bigg)
\nonumber\\
    & \qquad
        + \iint_{K\times\R} s^2\Phi_R(s) \, \chi(s|\OVZ) \,ds \,dx.
\label{E:UJN}
\end{align}
Since $\chi(s|\OVZ)\in L^1(K\times\R)$ there exists $R>0$ such that
$$
    \iint_{K\times\R} s^2\Phi_R(s) \, \chi(s|\OVZ) \,ds \,dx
        \LS \EPS/3.
$$
Moreover, we can find $N_1>0$ such that
$$
    \sup_{n\GS N_1} \bigg|
        \iint_{K\times\R} s^2 \chi(s|\OVZ^n) \,ds \,dx
            -\iint_{K\times\R} s^2 \chi(s|\OVZ) \,ds \,dx \bigg|
        \LS \EPS/3,
$$
by assumption \eqref{E:ENERCONV} of convergence of the initial
total
energies. For the remaining term on the right-hand side of
\eqref{E:UJN}, we define the function
$$
    \eta_R(\BA)
         := \int_\R s^2\varphi_R(s)\, \chi(s|\BA) \,ds
    \quad\text{for $\BA\in H$,}
$$
which is continuous and can be estimated as in \eqref{E:ETAEST}.
Therefore
$$
    \eta_R(\BA)
        \LS C_R\Big( 1+\rho(\BA)^{2\theta\lambda} \Big)
    \quad\text{for all $\BA\in H$,}
$$
with $C_R>0$ some constant.

Note that $\gamma>1$ implies $2\theta\lambda<\gamma$, so the
sequence $(\eta_R(\OVZ^n))$ is equi-integrable because of
\eqref{E:ENTRINEQ}. Since $\OVZ^n \longrightarrow \OVZ$ in measure
by assumption \eqref{E:STRONG}, we have
$$
    \eta_R(\OVZ^n) \longrightarrow \eta_R(\OVZ)
    \quad\text{strongly in $L^1(K)$.}
$$
Therefore there exists a number $N_2>0$ with
$$
    \sup_{n\GS N_2} \bigg| \iint_{K\times\R} s^2\varphi_R(s)\,
                \chi(s|\OVZ^n) \,ds \,dx
            - \iint_{K\times\R} s^2\varphi_R(s) \,
                \chi(s|\OVZ) \,ds \,dx \bigg|
    \LS \EPS/3.
$$
Combining all estimates, we obtain \eqref{E:CLA} with $N :=
\max(N_1,N_2)$.
\qed

Since the finite energy approximations $(\rho^n,u^n)$ are
themselves
entropy solutions of the isentropic Euler equations, we can use the
kinetic formulation, which implies the existence of
nonnegative measures $\mu^n$ such that
\begin{align}
    & \partial_t \Big( \eta_\psi(\BZ^n) A^n \Big)
        + \partial_x \Big( q_\psi(\BZ^n) A^n \Big)
        + \Big( \big(\rho u\; \eta_{\psi,\rho}-q_\psi\big)(\BZ^n)
            \Big) (\partial_x A^n)
\nonumber\\
    & \quad
        = -\int_\R \psi''(s)A^n \, \mu^n(ds,\cdot)
    \quad\text{in $\D'\big([0,\infty)\times\Omega\big)$,}
\label{E:WSX2}
\end{align}
for all test functions $\psi\in\D(\R)$. We also have
$\eta_\psi(\BZ^n(0, \cdot)) = \eta_\psi(\OVZ^n)$ in the
distribution
sense. Since the measures $\mu^n$ are uniformly bounded:
\begin{align}
    \iint_{[0,\infty)\times\Omega} \int_\R  A^n(x)\,
            \mu^n(ds,dx,dt)
        & \LS \int_\R \Big( \HA\OVR^n (\OVU^n)^2
            + U(\OVR^n) \Big) A^n \,dx
\nonumber\\
        & \LS 2\OVE
    \quad\text{for all $n$}
\label{E:BDMEASN}
\end{align}
(see \eqref{E:ENERCONV}), we obtain that along a suitable
subsequence (still denoted by $\mu^n$)
$$
    A^n\mu^n \WEAK A\mu
    \quad\text{weak-$\star$ in $M\Big( \big(
        [0,\infty)\times\bar{\Omega} \big) \times \R\Big)$.}
$$
Recall that $A^n$ converges uniformly to $A$, by construction. After
extracting another subsequence if necessary, we may also assume that
the sequence $(\rho^n,u^n)$ generates a Young measure $\nu$ as
introduced in Proposition~\ref{P:FIRST}. Using Lemmas~\ref{L:PSIIND}
\& \ref{L:CONID}, we can then pass to the limit in equation
\eqref{E:WSX2} and obtain
\begin{align}
    \partial_t \Big(\LA\eta_\psi\RA A\Big)
        + \partial_x \Big(\LA q_\psi\RA A\Big)
        + \big\LA \rho u\; \eta_{\psi,\rho}-q_\psi\big\RA
            (\partial_x A)
        &= -\int_\R \psi''(s)A \,\mu(ds,\cdot),
\nonumber\\
    \LA\eta_\psi\RA(0,\cdot)
        &= \eta_\psi(\OVZ)
\label{E:KINN}
\end{align}
in $\D'([0,\infty)\times \Omega)$ for all test functions
$\psi\in\D(\R)$. In this sense, the Young measure $\nu$ is a
measure-valued solution of the isentropic Euler equations
\eqref{E:CONS}. In the next subsection we are going to show that
\eqref{E:KINN} extends to weight functions $\psi$ that have
subquadratic growth at infinity. This will in particular imply that
the initial data $(\OVR, \OVU)$ is attained in the distribution
sense.


\subsection{Equi-integrability of the energy}\label{SS:EQUII}

Here is an extension of Lemma~\ref{L:PSIIND}.

\begin{proposition}[Higher integrability of the
energy]\label{P:NONL} Assume that the sequence $(\rho^n,u^n)$ of
finite energy approximations constructed in Section~\ref{SS:APPROX}
generates a Young measure $\nu$ as explained in
Proposition~\ref{P:FIRST}. Consider the sequence $(\BZ^n)$ of
Riemann invariants associated with $(\rho^n,u^n)$. For any weight
$\psi\in C^2(\R)$ with subcubic growth at infinity, we then obtain
\begin{equation}
    \eta_\psi(\BZ^n)A^n \WEAK \LA\eta_\psi\RA A
    \quad\text{weakly in
        $L^1_\LOC\big([0,\infty)\times\bar{\Omega}\big)$.}
\label{E:CONSQU1}
\end{equation}
Moreover, if $\psi$ has subquadratic growth at infinity, then
\begin{equation}
    \begin{aligned}
        q_\psi(\BZ^n)A^n &\WEAK \LA q_\psi\RA A
\\
        (\rho u\; \eta_{\psi,\rho})(\BZ^n)A^n
            &\WEAK \LA \rho u\; \eta_{\psi,\rho} \RA A
    \end{aligned}
    \quad\text{weakly in
        $L^1_\LOC\big([0,\infty)\times\bar{\Omega}\big)$.}
\label{E:CONSQU2}
\end{equation}
\end{proposition}

Proposition~\ref{P:NONL} shows that in \eqref{E:KINN} we can allow
weight functions $\psi$ that do not have compact support, but grow
subquadratically at infinity. In particular, we can choose
$\psi(s)=1$ or $\psi(s)=s$, and obtain the analogue of the
continuity and momentum equation in \eqref{E:CONS} for the
measure-valued solution $\nu$.

The following lemma is a generalization of results from \cite{LPT,
Mi}.

\begin{lemma}
\label{L:BFLUX}
Let $(\rho^n,u^n)$ be the sequence of finite energy
approximations from Section~\ref{SS:APPROX}. Then there exists a
constant $C>0$ such that for all $T>0$
\begin{equation}
    \sup_n\; \ESUP_{y\in\Omega} \Bigg\{ A^n(y)
        \int_{[0,T]} \Big( \rho^n |u^n|^3 +
(\rho^n)^{\gamma+\theta}
            \Big)(t,y) \,dt \Bigg\} \LS C.
\label{E:FLUXBD}
\end{equation}
\end{lemma}

{\em Proof.}
As explained at the beginning of Section~\ref{SS:MEASURE}, for any
$n$ there exists a nonnegative measure $\mu^n$ such that in the
distribution sense
\begin{align}
    & \partial_t \Big( \chi(\cdot|\rho^n,u^n) A^n \Big)
        + \partial_x \Big( \sigma(\cdot|\rho^n,u^n) A^n \Big)
\nonumber\\
    & \quad
        + \Big( \rho^n u^n\; \chi_{,\rho}(\cdot|\rho^n,u^n)
            -\sigma(\cdot|\rho^n,u^n) \Big) (\partial_x A^n)
        = -\partial^2_s(A^n\mu^n).
\label{E:KINETICN}
\end{align}
We now integrate \eqref{E:KINETICN} against the function
$$
    \IND_{[0,T]\times[y,\infty)}(t,x) \psi(s)
$$
with $\psi(s) := \HA s|s|$ for $s\in\R$. Using a standard
approximation argument, we obtain that
for almost every $T\in[0,\infty)$ and $y\in\Omega$
\begin{align}
    & A^n(y) \int_{[0,T]} q_\psi(\rho^n,u^n)(t,y) \,dt
\nonumber\\
    & \quad
        = \int_{[y,\infty)} \eta_\psi(\rho^n,u^n)(T,x)
            A^n(x) \,dx
        - \int_{[y,\infty)} \eta_\psi(\rho^n,u^n)(0,x)
            A^n(x) \,dx
\nonumber\\
    & \qquad
        + \iint_{[0,T]\times[y,\infty)}
            \Big( \rho^n u^n\; \eta_{\psi,\rho}(\rho^n,u^n)
                -q_\psi(\rho^n,u^n) \Big)(t,x)
                    \big(\partial_x A^n\big)(x) \,dx \,dt
\nonumber\\
    & \qquad
        + \iint_{[0,T]\times[y,\infty)}
            \SIGN(s) A^n(x) \,\mu^n(ds,dx,dt).
\label{E:DSXC}
\end{align}
As usual, the entropy/entropy-flux pair $(\eta_\psi, q_\psi)$ is
defined by \eqref{E:ETAPSI}. Now
$$
    \bigg| \iint_{[0,T]\times[y,\infty)}
        \SIGN(s) A^n(x) \,\mu^n(ds,dx,dt) \bigg|
            \LS 2\OVE
$$
for all $n$ because of \eqref{E:BDMEASN}. Moreover, since for all
finite energy approximations the total energy is nonincreasing in
time, we can estimate for $t\in\{0,T\}$
\begin{align*}
    \bigg| \int_{[y,\infty)} \eta_\psi(\rho^n,u^n)(t,x)
            A^n(x) \,dx \bigg|
    & \LS \int_\Omega \Big( \HA\rho^n(u^n)^2 + U(\rho^n)
        \Big)(t,x) A^n(x) \,dx
\\
    & \LS \int_\Omega \Big( \HA\OVR^n(\OVU^n)^2 + U(\OVR^n)
        \Big)(x) A^n(x) \,dx,
\end{align*}
which for all $n$ is bounded by $2\OVE$ (see \eqref{E:ENTRINEQ} and
\eqref{E:ENERCONV}). Recall that the total energy is the second
$s$-moment of the entropy kernel. For the third integral on the
right-hand side of \eqref{E:DSXC}, a computation based on
\eqref{E:TED} yields
$$
    \rho^n u^n\; \eta_{\psi,\rho}(\rho^n,u^n)-q_\psi(\rho^n,u^n)
        = -\theta (\rho^n)^{\gamma+\theta}
            \frac{\big(1-u^n/(\rho^n)^\theta\big)^{\lambda+2}_+}
                {(\lambda+1)(\lambda+2)}.
$$
This quantity is nonpositive and bounded below by $-C(\rho^n)
^{\gamma+\theta}$, with $C>0$ some constant. Finally, we use the
fact that there exists $\delta>0$ such that
$$
    q_\psi(\rho^n,u^n) \GS \delta \Big( \rho^n |u^n|^3
        + (\rho^n)^{\gamma+\theta} \Big)
    \quad\text{for all $(\rho^n,u^n)$.}
$$
We refer the reader to \cite{LPT} for a proof. Combining all
estimates, we find
\begin{equation}
    Q^n(y) \LS \frac{6\OVE}{\delta}
        + \frac{C}{\delta} \int_{[y,\infty)}
            \frac{\big(\partial_x A^n(x)\big)_-}{A^n(x)}
                \, Q^n(x) \,dx
\label{E:LGB}
\end{equation}
for almost all $y\in\Omega$, where
$$
    Q^n(y) := A^n(y) \int_{[0,T]} \Big( \rho^n |u^n|^3
        + (\rho^n)^{\gamma+\theta} \Big)(t,y) \,dt.
$$
Note that for every $n$, the functions $(\rho^n,u^n)$ and
$Q^n$ are compactly supported, so the integral in \eqref{E:LGB} is
well-defined. Then Gronwall's lemma implies
\begin{equation}
    Q^n(y) \LS \frac{6\OVE}{\delta} \exp\!\Bigg(
        \frac{C}{\delta} \int_{[y,\infty)}
            \frac{\big(\partial_x A^n(x)\big)_-}{A^n(x)}\,dx
                \Bigg)
    \quad\text{for a.e.\ $y\in\Omega$.}
\label{E:HGML}
\end{equation}
For nozzle flows, the right-hand side of \eqref{E:HGML} can
be bounded independently of $y$ and $n$, by assumption
\eqref{E:GEON} and the choice of $A^n$. For spherically symmetric
flows, the weight $A^n$ is strictly increasing, so the integral in
\eqref{E:HGML} vanishes.
\qed

{\em Proof of Proposition~\ref{P:NONL}.}
Let $p:=(\gamma+\theta)/\gamma$ such that $p>1$. Then
\begin{equation}
    \sup_n \iint_{[0,T]\times K} \Big( \rho^n(u^n)^2
        +(\rho^n)^\gamma \Big)^p A^n \,dx \,dt
    \LS C
\label{E:HITE}
\end{equation}
for all $T>0$ and $K\subset\bar{\Omega}$ compact, with $C>0$ some
constant: Note first that
\begin{align}
    & A^n \int_{[0,T]} \Big( \rho^n (u^n)^2 \Big)^p \,dt
\nonumber\\
    & \quad
        \LS \bigg( A^n \int_{[0,T]} \rho^n |u^n|^3 \,dt
            \bigg)^{\!(3\gamma-1)/3\gamma}
        \bigg( A^n \int_{[0,T]} (\rho^n)^{\gamma+\theta}
            \,dt \bigg)^{\!1/3\gamma},
\label{E:LOKH}
\end{align}
by H\"{o}lder inequality. For the internal energy, we have the
trivial identity
\begin{equation}
    A^n \int_{[0,T]} \Big( (\rho^n)^\gamma \Big)^p \,dt
        = A^n \int_{[0,T]} (\rho^n)^{\gamma+\theta} \,dt.
\label{E:LOKH2}
\end{equation}
Since the right-hand sides of both \eqref{E:LOKH} and
\eqref{E:LOKH2} are bounded independently of $x$ and $n$ because of
Lemma~\ref{L:BFLUX}, the bound \eqref{E:HITE} follows immediately
after integrating over $K$. Similarly, we can use the H\"{o}lder
inequality to prove
\begin{equation}
    \sup_n \iint_{[0,T]\times K} (\rho^n)^\gamma |u^n|
        A^n \,dx \,dt \LS C
\label{E:LECK}
\end{equation}
for some constant $C>0$. Indeed, we have
$$
    A^n \int_{[0,T]} (\rho^n)^\gamma |u^n| \,dt
        \LS \bigg( A^n \int_{[0,T]} (\rho^n)^{\gamma+\theta}
            \,dt \bigg)^{\!2/3}
        \bigg( A^n \int_{[0,T]} \rho^n|u^n|^3 \,dt
            \bigg)^{\!1/3},
$$
which is bounded uniformly. Integrating over $K$, we
obtain \eqref{E:LECK}. Thus
\begin{equation}
    \sup_n \iint_{[0,T]\times K} \bigg( \int_\R s^2\chi(s|\rho^n,
        u^n) \,ds \bigg) |u^n| A^n \,dx \,dt \LS C
\label{E:EDI1}
\end{equation}
because the second $s$-moment of $\chi$ is given by the total
energy.

Let again $\psi(s):=s|s|$ for $s\in\R$. Then formulas
\eqref{E:CHISIG} \& \eqref{E:ENTROPY} imply
$$
    \theta \int_\R |s|^3 \chi(s|\rho^n,u^n) \,ds
        = q_\psi(\rho^n,u^n)
        - (1-\theta) u^n \int_\R s|s|
            \chi(s|\rho^n,u^n) \,dt.
$$
The first term on the right-hand side can be controlled using the
argument of Lemma~\ref{L:BFLUX} (see \eqref{E:DSXC}). For the
second
term, we can use \eqref{E:EDI1}. This yields
\begin{equation}
    \sup_n \iint_{[0,T]\times K} \bigg( \int_\R |s|^3\chi(s|\rho^n,
        u^n) \,ds \bigg) A^n \,dx \,dt \LS C,
\label{E:EDI2}
\end{equation}
with $C>0$ some constant. Combining \eqref{E:EDI1} \&
\eqref{E:EDI2}, we obtain the convergence of
$\eta_\psi(\BZ^n)$ and $q_\psi(\BZ^n)$ for unbounded $\psi$ by
standard arguments.

To prove the last statement in \eqref{E:CONSQU2}, note that
\begin{align}
    \rho^n u^n\; \eta_{\psi,\rho}(\rho^n,u^n)
        & = u^n \int_\R \psi(s) \chi(s|\rho^n,u^n) \,ds
\nonumber\\
        & \quad
            + \theta u^n \int_\R \psi'(s)(s-u^n) \chi(s|\rho^n,
                u^n) \,ds.
\label{E:TED}
\end{align}
Using \eqref{E:FLUXBD} and \eqref{E:EDI1}, we can control the
right-hand side of \eqref{E:TED} uniformly in $n$, for all $\psi$
with at most quadratic growth. This completes the proof.
\qed


\subsection{Compensated compactness}\label{SS:COMPC}

We have the following crucial result.

\begin{lemma}[$\DIV$-$\CURL$-commutator]\label{L:COM}
Assume that the sequence $(\rho^n,u^n)$ of finite energy
approximations constructed in Section~\ref{SS:APPROX} generates a
Young measure $\nu$. Then almost everywhere in $[0,\infty)
\times\Omega$ we have
\begin{align*}
    \LA \chi(s) \sigma(s') - \sigma(s) \chi(s') \RA
            - \LA\chi(s)\RA \LA\sigma(s')\RA
        & + \LA\sigma(s)\RA \LA\chi(s')\RA = 0
\\
        & \qquad \text{for a.e.\ $(s,s')\in\R^2$.}
\end{align*}
\end{lemma}

{\em Proof.}
For any test functions $\psi,\psi'\in\D(\R)$ define the
entropy/entropy-flux pairs $(\eta_\psi, q_\psi)$ and
$(\eta_{\psi'},
q_{\psi'})$ as in \eqref{E:ETAPSI}. According to
Lemma~\ref{L:PSIIND} we have
\begin{align}
    \left.\begin{aligned}
        \eta_\psi(\BZ^n) &\WEAK \LA\eta_\psi\RA
\\
        q_\psi(\BZ^n) &\WEAK \LA q_\psi\RA
    \end{aligned}\right.
    & \quad\text{weakly in $L^{\gamma+1}_\LOC
        \big([0,\infty)\times\Omega\big)$,}
\label{E:WSX1}
\end{align}
as well as
\begin{equation}
    (\rho u\; \eta_{\psi,\rho})(\BZ^n)
        \WEAK \LA \rho u\; \eta_{\psi,\rho} \RA
    \quad\text{weakly in $L^2_\LOC
        \big([0,\infty)\times\Omega\big)$.}
\label{E:WSX1A}
\end{equation}
The same convergence holds for the pair $(\eta_{\psi'},
q_{\psi'})$. Moreover, we have
\begin{align}
    \left.\begin{aligned}
        \eta_\psi(\BZ^n) q_{\psi'}(\BZ^n)
            &\WEAK \LA\eta_\psi q_{\psi'}\RA
\\
        q_\psi(\BZ^n) \eta_{\psi'}(\BZ^n)
            &\WEAK \LA q_\psi \eta_{\psi'}\RA
    \end{aligned}\right.
    & \quad\text{weakly in
$L^1_\LOC\big([0,\infty)\times\Omega\big)$.}
\label{E:WSX3}
\end{align}
Recall that for all $\psi\in\D(\R)$, the sequence $(\BZ^n)$
satisfies
\begin{align}
    & \partial_t \Big( \eta_\psi(\BZ^n) A^n \Big)
        + \partial_x \Big( q_\psi(\BZ^n) A^n \Big)
        + \Big( \big(\rho u\; \eta_{\psi,\rho}-q_\psi\big)(\BZ^n)
            \Big) (\partial_x A^n)
\nonumber\\
    & \quad
        = -\int_\R \psi''(s)A^n \, \mu^n(ds,\cdot)
    \quad\text{in $\D'\big([0,\infty)\times\Omega\big)$.}
\label{E:WSX2N}
\end{align}
By \eqref{E:BDMEASN}, the right-hand side of \eqref{E:WSX2N} is
bounded in $\MEAS([0,\infty)\times \Omega)$. Moreover, by
\eqref{E:WSX1} \& \eqref{E:WSX1A} and the divergence form of the
left-hand side of \eqref{E:WSX2N}:
\begin{align*}
    & \bigg( \int_\R \psi''(s)A^n \,\mu^n(ds,\cdot) \bigg)
\\
    & \qquad
    \begin{array}{l}
        \text{is pre-compact in
            $W_\LOC^{-1,r}\big([0,\infty)\times
            \Omega\big)$ for $1\LS r<2$}
\\[-0.5ex]
        \text{and uniformly bounded in $W^{-1,\gamma+1}
            _\LOC\big([0,\infty)\times\Omega\big)$.}
    \end{array}
\end{align*}
We used Sobolev embedding. Since $\gamma+1>2$, Murat's
Lemma \cite{Murat2} yields
$$
    \left(\int_\R \psi''(s)A^n \,\mu^n(ds,\cdot)\right)
    \quad\text{is pre-compact in $H^{-1}_\LOC\big([0,\infty)\times
        \Omega\big)$.}
$$
The same arguments apply to the entropy/entropy-flux pair
$(\eta_{\psi'},q_{\psi'})$.

We now use the $\DIV$-$\CURL$-Lemma (see \cite{Murat, Tartar}),
which gives the identity
\begin{equation}
    \LA-\eta_\psi q_{\psi'} + q_\psi \eta_{\psi'}\RA
        + \LA\eta_\psi\RA \LA q_{\psi'}\RA
        - \LA q_\psi\RA \LA\eta_{\psi'}\RA = 0
    \quad\text{in $\D'\big([0,\infty)\times\Omega\big)$.}
\label{E:IKJ}
\end{equation}
By \eqref{E:WSX1} and \eqref{E:WSX3}, the commutator is in
$L^1_\LOC([0,\infty)\times\Omega)$, so \eqref{E:IKJ} holds
pointwise almost everywhere. On the other hand, by \eqref{E:ETAPSI}
we have
\begin{align*}
    & \LA-\eta_\psi q_{\psi'} + q_\psi \eta_{\psi'}\RA
        + \LA\eta_\psi\RA\LA q_{\psi'}\RA
        - \LA q_\psi\RA\LA\eta_{\psi'}\RA
\\
    & \quad
        = \iint_{\R^2}
            \Big( \LA-\chi(s)\sigma(s')+\sigma(s)\chi(s')\RA
                + \LA\chi(s)\RA \LA\sigma(s')\RA
                - \LA\sigma(s)\RA \LA\chi(s')\RA \Big)
\\
    & \quad\hspace{20em}
        \psi(s) \psi'(s') \,ds \,ds'.
\end{align*}
Since $\psi,\psi'$ were arbitrary, the integrand must vanish for
almost all $(s,s')$.
\qed


\section{Strong convergence and finite energy
solutions}\label{S:STRONG}

In the previous section, we showed that a subsequence of the finite
energy approximate solutions $(\rho^n,u^n)$ converges to a
measure-valued solution of the isentropic Euler equations. In this
section, we improve this result by showing that the Young measure
constructed in Proposition~\ref{P:FIRST} is concentrated for a.e.\
$(t,x)\in[0,\infty) \times\Omega$. This implies the existence of
measurable functions $(\rho,u)$, which form a weak solution in the
sense of Definition~\ref{D:SOLUTION}.

\subsection{Reduction of the Young measure}\label{SS:LEMMAS}

We first introduce some notation.

\begin{definition}\label{D:ADM}
Consider $\nu\in\PROB(\CC)$ such that $\LA W\RA$ is
finite, where
$$
    \LA\varphi\RA := \int_\CC \varphi(\BA) \,\nu(d\BA)
$$
 for all $\varphi := \varphi_\BD W$ with $\varphi_\BD \in
C_\BD(\CC)$. The measure $\nu$ is called an entropy admissible Young
measure if for almost every $(s,s')\in\R^2$ we have
\begin{equation}
    \LA \chi(s) \sigma(s') - \sigma(s) \chi(s') \RA
            -\LA\chi(s)\RA \LA\sigma(s')\RA
            +\LA\sigma(s)\RA \LA\chi(s')\RA
        = 0.
\label{E:COMR}
\end{equation}
\end{definition}

Entropy admissible measures have a very particular structure:

\begin{theorem}[Reduction of Young measures]
\label{T:CONC} If $\nu$ is an entropy admissible Young measure,
then
the support of $\nu$ is either a single point of\, $H$ or a subset
of the vacuum line $V$.
\end{theorem}

As shown in Proposition~\ref{P:FIRST} and Lemma~\ref{L:COM}, the
sequence $(\rho^n,u^n)$ of finite energy approximate solutions
constructed in Subsection~\ref{SS:APPROX}, generates a Young measure
with the property that for almost every $(t,x) \in [0,
\infty)\times\Omega$ the measure $\nu_{(t,x)}$ is entropy admissible
in the sense of Definition~\ref{D:ADM}. We can therefore apply
Theorem~\ref{T:CONC} at each point: For all $(t,x)$ where
$\nu_{(t,x)}$ is not supported in the vacuum, we have $\nu_{(t,x)} =
\delta_{\BZ(t,x)}$ for some $\BZ(t,x)\in H$, thus
\begin{align}
    \LA\eta_\psi\RA(t,x) &= \eta_\psi\big(\BZ(t,x)\big),
\nonumber\\
    \LA q_\psi\RA(t,x) &= q_\psi\big(\BZ(t,x)\big),
\label{E:ENDL}\\
    \LA\rho u\; \eta_{\psi,\rho}-q_\psi\RA(t,x)
        &= \big( \rho u\; \eta_{\psi,\rho}-q_\psi\big)
            \big(\BZ(t,x)\big)
\nonumber
\end{align}
for all admissible weight functions $\psi$. If $\nu_{(t,x)}$ is
supported in $V$, then
$$
    \LA\eta_\psi\RA(t,x)
        = \LA q_\psi\RA(t,x)
        = \LA\rho u\; \eta_{\psi,\rho}-q_\psi\RA(t,x)
        = 0
$$
since the integrands vanish in the vacuum, see \eqref{E:ERT} and
\eqref{E:QEST}. For those points we define $\BZ(t,x) := (0,0)$ and
obtain again \eqref{E:ENDL}. The Young measure $\nu$ is a
measure-valued solution of the isentropic Euler equations in the
sense \eqref{E:KINN}. With $\BZ\colon[0,\infty) \times \Omega
\longrightarrow\CC$ defined above \eqref{E:KINN} takes the form
\begin{align}
    \partial_t \Big(\eta_\psi(\BZ)A\Big)
        + \partial_x \Big(q_\psi(\BZ)A\Big)
        + \Big( \big(\rho u\; \eta_{\psi,\rho}
            & -q_\psi\big)(\BZ) \Big) (\partial_x A)
\nonumber\\
        &= -\int_\R \psi''(s)A \,\mu(ds,\cdot),
\nonumber\\
    \eta_\psi\big(\BZ(0,\cdot)\big)
        &= \eta_\psi(\OVZ)
\label{E:KINN2}
\end{align}
in $\D'([0,\infty)\times \Omega)$ for all admissible weight
functions $\psi$.

Consider now the functions $(\rho,u)$ that are related to $\BZ$ via
\eqref{E:RIEMANN}. Then \eqref{E:KINN2} shows that $(\rho,u)$ is an
entropy solution in the sense of Definition~\ref{D:SOLUTION}, which
proves our main Theorem~\ref{T:MAIN}. Observe that in
Proposition~\ref{P:NONL} we can allow functions $\psi$ with
quadratic growth in the entropy $\LA\eta_\psi\RA$, but only
subquadratic growth is acceptable for the entropy-flux $\LA
q_\psi\RA$. Since for the finite energy approximate solutions the
total energy is nonincreasing in time, the same is true for the
limit functions $(\rho,u)$/ We therefore have
$$
    \int_\Omega \big( \HA\rho u^2+U(\rho) \big)(t_2,x) \,dx
        \LS \int_\Omega \big( \HA\rho u^2+U(\rho) \big)(t_1,x) \,dx
$$
for almost every $t_2\GS t_1$. Note, however, that while the
argument of Lemma~\ref{L:BFLUX} can be used to derive a uniform
$L^1$-bound for the total energy fluxes
$$
    \big( \HA\rho^n (u^n)^2+Q(\rho^n) \big) u^n A^n,
$$
we cannot prove that their limit is given by
$$
    \big( \HA\rho u^2+Q(\rho) \big) u A
$$
since concentrations might occur. As a consequence, we do not
know whether the local energy balance (that is, \eqref{E:CENER}
with an inequality) is satisfied.

The rest of this section is devoted to the  proof of
Theorem~\ref{T:CONC}.

\begin{lemma}
\label{L:CHICON}
Given an entropy admissible Young measure $\nu$,
consider the map $s \in \R \mapsto \LA \chi(s) \RA$. Then,
$\LA\chi\RA\in C^\alpha(\R)$ for all $\alpha\in[0,\lambda]$, and so
the set
$$
    \SC := \big\{ s\in\R \colon \LA\chi(s)\RA > 0 \big\}
$$
is open. If\, $\SC$ is empty, then $\nu(H)=0$. If\, $\SC$ is
nonempty, define numbers $\UZ:=\inf\SC$ and $\OZ:=\sup\SC$ (both
possibly unbounded). Then $\SC=(\UZ,\OZ)$ and
\begin{equation}
    \SPT\nu\cap\big\{ \BA\in H\colon
        \text{$\UA<\UZ$ or $\OZ<\OA$} \big\} = 0.
\label{E:TRIANGLE}
\end{equation}
\end{lemma}

{\em Proof.}
Note that the function $f(t) := (1-t^2)^\lambda_+$ is bounded
and H\"{o}lder continuous with H\"{o}lder exponent $\lambda$. We
write the entropy kernel in the form
\begin{equation}
    \chi(s|\BA) = \rho(\BA)^{2\theta\lambda}\,
        f\bigg( \frac{s-u(\BA)}{\rho(\BA)^\theta} \bigg)
    \qquad\text{for $(s,\BA)\in\R\times\CC$,}
\label{E:SCALING}
\end{equation}
where $\rho(\BA)$ and $u(\BA)$ are defined by \eqref{E:RIEMANN}. We
then obtain
\begin{align*}
    \sup_{s\neq s'} \frac{|\chi(s|\BA)-\chi(s'|\BA)|}
            {|s-s'|^\alpha}
        & = \rho(\BA)^{(2\lambda-\alpha)\theta}\;
            \sup_{t\neq t'} \frac{|f(t)-f(t')|}{|t-t'|^\alpha}
\\
        & \LS C \rho(\BA)^{(2\lambda-\alpha)\theta}\;,
\end{align*}
with $C>0$ some constant that does not depend on $\BA$. We also
have
$$
    \sup_{s\in\R} |\chi(s|\BA) | \LS \rho(\BA)^{2\lambda\theta}.
$$

Since $0<(2\lambda-\alpha)\theta<1$ for all $\alpha\in[0,\lambda]$,
we can now estimate
\begin{align*}
    \sup_{s\neq s'} \frac{ \big| \LA\chi(s)\RA-\LA\chi(s')\RA
\big|}
            {|s-s'|^\alpha}
        & = \sup_{s\neq s'} |s-s'|^{-\alpha} \, { \DST \bigg|
\int_\CC
\chi(s|\BA)
            \,\nu(d\BA) - \int_\CC \chi(s'|\BA) \,\nu(d\BA) \bigg|}
   \\
    & \LS \int_\CC \sup_{s\neq s'}
\frac{|\chi(s|\BA)-\chi(s'|\BA)|}
            {|s-s'|^\alpha} \,\nu(d\BA)
   \\
   & \LS C \int_\CC W(\BA) \,\nu(d\BA),
\end{align*}
which is finite by assumption on $\nu$. The function $\LA\chi\RA$
is bounded:
\begin{align*}
    \sup_{s\in\R} |\LA\chi(s)\RA|
        & = \sup_{s\in\R} \bigg| \int_\CC \chi(s|\BA)
            \,\nu(d\BA) \bigg|
\\
        & \LS \int_\CC \sup_{s\in\R} |\chi(s|\BA)| \,\nu(d\BA)
        \LS \int_\CC W(\BA) \,\nu(d\BA).
\end{align*}
This shows that $\LA\chi\RA\in C^\alpha(\R)$ for all $\alpha\in[0,
\lambda]$, so $\SC$ is well-defined and open.

We show next that $\SC$ can be represented in the form
\begin{equation}
    \SC = \bigcup_{\BA\in\SPT\nu\cap H} (\UA,\OA).
\label{E:SP}
\end{equation}
Indeed assume that $\BA\in\SPT\nu\cap H$. Then we have
$\nu(B_r(\BA)
\cap H)>0$ for all $r>0$, by definition of support of a measure.
Therefore we obtain
$$
    \LA\chi(s)\RA
        \GS \int_{B_r(\BA)} \chi(s|\BA') \,d\nu(\BA') >0
$$
at least for all $s\in\R$ with the property that $\chi(s|\BA')>0$
for all $\BA'\in B_r(\BA)$. This implies $(\UA+r,\OA-r) \subset
\SC$. Since $r>0$ and $\BA$ were arbitrary, we get the $\supset$
inclusion in \eqref{E:SP}. For the converse direction, suppose that
\begin{equation}
    \LA\chi(s)\RA
        = \int_\CC \chi(s|\BA') \,d\nu(\BA') > 0
\label{E:AS}
\end{equation}
for some $s\in\R$. Since $\chi$ vanishes in the vacuum, in
\eqref{E:AS} we can restrict integration to $H$. Then $\nu(\{\BA\in
H\colon \UA<s<\OA \}) > 0$, so there exists at least one point
$\BA\in\SPT\nu$ in that set. Then $s\in (\UA,\OA)$, and \eqref{E:SP}
follows. If now $\SC$ is empty, then \eqref{E:SP} implies that
$\SPT\nu\cap H=\varnothing$, thus $\nu(H)=0$.

Let us now assume that $\SC$ is nonempty. We define $\UZ,\OZ$
as in the statement of the lemma. Then we argue by contradiction
and assume that $\SC$ is disconnected. Since $\SC$ is open, there
exist numbers $\UZ < \UC \LS \OC < \OZ$ and $\EPS>0$ such that
$$
\begin{cases}
    \;\LA\chi(s)\RA=0
        & \text{for $s\in[\UC,\OC]$,}
\\
    \;\LA\chi(s)\RA>0
        & \text{for $s\in (\UC-\EPS,\UC)
            \cup (\OC,\OC+\EPS)$.}
\end{cases}
$$
In view of \eqref{E:SP}, this implies that
\begin{equation}
    \SPT\nu \cap \big\{ \BA\in H\colon \text{$\UC<\OA$ and
            $\UA<\OC$} \big\}
        = \varnothing.
\label{E:ESD}
\end{equation}
Choosing $s\in (\UC-\EPS,\UC)$ and $s'\in (\OC,\OC+\EPS)$ we use
assumption \eqref{E:COMR} in the form
\begin{equation}
    \frac{\LA-\chi(s)\sigma(s')+\sigma(s)\chi(s')\RA}
            {\LA\chi(s)\RA \LA\chi(s')\RA}
        = \frac{\LA\sigma(s')\RA}{\LA\chi(s')\RA}
            -\frac{\LA\sigma(s)\RA}{\LA\chi(s)\RA},
\label{E:TRF}
\end{equation}
which is well-defined since $\LA\chi(s)\RA \LA\chi(s')\RA>0$. Now
note that $\chi(s|\BA)\chi(s'|\BA)=0$ for all $\BA\in\SPT\nu$, by
\eqref{E:ESD}.
%
%
We obtain
$$
    -\chi(s|\BA)\sigma(s'|\BA)
        +\sigma(s|\BA)\chi(s'|\BA) = 0
    \quad\text{for all $\BA\in\SPT\nu$,}
$$
so the left-hand side of \eqref{E:TRF} vanishes. For the right-hand
side we can estimate
$$
    \frac{\LA\sigma(s)\RA}{\LA\chi(s)\RA}
        = \theta s \frac{\LA\chi(s)\RA}{\LA\chi(s)\RA}
            + (1-\theta) \frac{\LA u\chi(s)\RA}{\LA\chi(s)\RA}
        \LS \theta s+(1-\theta)\UC
        < \UC.
$$
Here, we have used that on the one hand
$$
    \SPT\chi(s|\cdot)\cap\SPT\nu
        \subset \big\{ \BA\in H\colon \OA\LS\UC \big\}\cup V
        \subset \big\{ \BA\in H\colon u(\BA)\LS \UC \big\}
            \cup V
$$
in view of \eqref{E:ESD}
and, on
the other hand, $\nu$ can not be entirely concentrated at one point
where $\chi(s|\BA)=0$ since $\LA\chi(s)\RA>0$.

With the analogous estimate
$$
    \frac{\LA\sigma(s')\RA}{\LA\chi(s')\RA}
        = \theta s' \frac{\LA\chi(s')\RA}{\LA\chi(s')\RA}
            + (1-\theta) \frac{\LA u\chi(s')\RA}{\LA\chi(s')\RA}
        \GS \theta s'+(1-\theta)\OC
        > \OC,
$$
we obtain from \eqref{E:TRF} that $0 > \OC-\UC \GS 0$, which is a
contradiction.
\qed

\subsection{Expansion of the entropy kernels}
\label{S:KERN}

In order to establish that the probability measure of
Theorem~~\ref{T:CONC} is concentrated at one point, we must
understand how the entropy kernels behave under fractional
differentiation with respect to $s$. For $\lambda>0$ and suitable
functions $f\colon \R\longrightarrow\R$ we define the operators
\begin{equation}
    \DC f
        := \F^{-1}\big(|\cdot|^{\lambda+1} \F f\big),
        \qquad
    \DS f
        := \F^{-1}\big(i|\cdot|^\lambda \SIGN(\cdot) \F f\big)
\label{E:FD}
\end{equation}
in distributional sense, where $\F$ denotes the Fourier transform.
We have
\begin{align}
    \DC f(s)
        &= \frac{d}{ds} \big( \DS f(s) \big),
\label{E:RFV1}\\
    \vphantom{\frac{d}{ds}}
    \DC \big( sf(s) \big)
        &= s\DC f(s) + (\lambda+1)\DS f(s).
\label{E:RFV2}
\end{align}
We now apply these operators to the function $f(s) :=
(1-s^2)^\lambda_+$ with $s\in\R$. According to \cite{GFS}, its
Fourier transform is given by
\begin{equation}
    \F f(z) :=
        2^\lambda\Gamma(\lambda+1) |z|^{-\lambda-1/2}
                J_{\lambda+1/2}(|z|)
\label{E:FOURIER}
\end{equation}
for all $z\in\R$, where $\Gamma$ denotes the Gamma function and
$J_{\lambda+1/2}$ is the Bessel function. Note that despite of the
singular factor in \eqref{E:FOURIER}, the function $\F f$ is
bounded, due to the decaying properties of the Bessel function. We
have
\begin{equation}
    \DS f = c\F ^{-1}\Big( |\cdot|^{-1/2} \F g \Big),
    \label{E:WEF}
\end{equation}
where $c$ is some constant and the function $g$ is defined for all
$z\in\R$ by
$$
    \F g(z) := i\SIGN(z) J_{\lambda+1/2}\big(|z|\big).
$$
The inverse Fourier transform of $|\cdot|^{-1/2}$ induces a
fractional integration operator, called Riesz potential (see
\cite{Stein}). Therefore \eqref{E:WEF} is equivalent to
\begin{equation}
    \DS f(s) = C|\cdot|^{-1/2} \star g(s),
    \qquad s\in\R,
\label{E:FALT}
\end{equation}
with $C$ some new constant. Since $\F g$ is an odd function, we can
express the inverse Fourier transform in terms of the inverse
Fourier Sine transform and obtain the following explicit formula
(see \cite{GRAD}):
\IGNORE{ 
\begin{align*}
    g(s)
        & = \frac{1}{\sqrt{2\pi}} \int_{-\infty}^\infty
            \F g(z) \exp(-izs) \,dz
\\
        & = \frac{-i}{\sqrt{2\pi}} \int_{-\infty}^\infty
            \F g(z) \sin(zs) \,dz
    \quad\text{for all $s\in\R$}
\end{align*}
since the $\cos$-term does not contribute. We decompose the
integral and get
} 
\begin{align}
    g(s)
        & = \sqrt{\frac{2}{\pi}} \SIGN(s) \int_0^\infty
            J_{\lambda+1/2}(z) \sin\!\big(z|s|\big) \,dz
\nonumber\\
        & = \sqrt{\frac{2}{\pi}} \SIGN(s)
\begin{cases}
    \DST\frac{\sin\!\Big( (\lambda+\HA) \arcsin|s| \Big)}
            {\sqrt{1-s^2} \vphantom{\Big)}},
        & \text{$|s|<1$,}
\\[2em]
    \DST\frac{\cos\!\Big( (\lambda+\HA)\frac{\pi}{2} \Big)}
            {\sqrt{s^2-1}\Big( |s|+\sqrt{s^2-1}
\Big)^{\lambda+1/2}},
        \qquad & \text{$|s|>1$.}
\end{cases}
\label{E:GLAM}
\end{align}
Note that $g$ decays like $|s|^{-(\lambda +3/2)}$ as
$|s|\rightarrow\infty$ and diverges only like $|1-|s||^{-1/2}$ as
$|s|\rightarrow 1$. This implies $g\in L^p(\R)$ for all
$p\in[1,2)$. By the Hardy-Little\-wood-Sobolev theorem (see
\cite{Stein}), we then have $\DS f\in L^q(\R)$ for all
$q\in(2,\infty)$. The singular behavior of $\DS f$ and $\DC f$ is
decribed in the following proposition.

\begin{proposition}[Fractional derivatives]
\label{P:EXP} Let $f(s)=(1-s^2)^\lambda_+$ for $s\in\R$, and
define the fractional derivatives $\DC f$ and $\DS f$ by
\eqref{E:FD}. Then there exist constants $A_i$, $i=1\ldots 4$, and
functions $r$, $q\in W^{1,p}(\R)$ for $p\in[2,\infty)$, such that
in the distribution sense we have the following expansions:
\begin{align*}
    \DS f(s) &=
          A_1 \Big( H(s+1)+H(s-1) \Big)
        + A_2 \Big( \CI(s+1)-\CI(s-1) \Big)
        + r(s),
\\
    \DC f(s) &=
          A_1 \Big( \delta(s+1)+\delta(s-1) \Big)
        + A_2 \Big( \PV(s+1)-\PV(s-1) \Big)
\\
    & \quad
        + A_3 \Big( H(s+1)-H(s-1) \Big)
        + A_4 \Big( \CI(s+1)+\CI(s-1) \Big)
        + q(s).
\end{align*}
Here $\delta$ is the Dirac measure, $\PV$ is the principal value
distribution, and $H$ denotes the Heaviside function. The
function $\CI$ is the Cosine integral
\begin{equation}
    \CI(s) := -\int_{|s|}^\infty \frac{\cos t}{t} \,dt
        = C + \log|s| +\int_0^{|s|}
            \frac{\cos t-1}{t} \,dt,
    \quad s\in\R,
\label{E:COSINE}
\end{equation}
with $C>0$ some constant. For simplicity, we will treat the
distributions $\delta$ and $\PV$ as if they were functions. The
coefficients $A_1$ and $A_2$ are not both equal to zero. Moreover,
if $\gamma =(M+2)/M$ with $M\in\N$ odd, then $A_2=A_4=0$.
%
%
\end{proposition}

\begin{remark}
Note that by Sobolev embedding, the remainders are H\"{o}lder
continuous: We have $r,q\in C^\alpha(\R)$ for all exponents
$\alpha\in[0,1)$. In particular, the functions are bounded.
Moreover, we get $r,q\in W^{1,p}_\LOC(\R)$ for all $p\in[1,
\infty)$.
\end{remark}

This expansion has been proved in slightly different form in
\cite{LPS, CL}, starting from an asymptotic formula for the Fourier
transform of $\DC f$. The main difference is that in \cite{LPS} the
logarithm $\log|\cdot|$ is used in place of $\CI$, which is not
totally accurate since the Fourier transform of $\DC f$ is a
bounded function, while the Fourier transform of the logarithm has
a pole at the origin. Recall that $\CI(s)$ behaves like $-\log|s|$
as $|s|\rightarrow 0$ and decays like $|s|^{-1}$ at infinity. We
remark in passing that it is possible to prove
Proposition~\ref{P:EXP} starting from identities \eqref{E:FALT} and
\eqref{E:GLAM}, thereby avoiding the Fourier transform altogether.
But we will not pursue this option here.

Proposition~\ref{P:EXP} is used to find expansions for the entropy
kernel. Note that
$$
    \chi(s|\BA) = \rho(\BA)^{2\theta\lambda}\, f\bigg(
        \frac{s-u(\BA)}{\rho(\BA)^\theta} \bigg),
    \qquad (s,\BA)\in\R\times\CC.
$$
Therefore the chain rule implies the identities
\begin{align}
    & \DS\chi(s|\BA)
\nonumber\\
    & \quad
        = \rho(\BA)^{\theta\lambda}
            \bigg( A_1 \Big( H(s-\UA)+H(s-\OA) \Big)
                + A_2 \Big( \CI(s-\UA)-\CI(s-\OA) \Big) \bigg)
\nonumber\\
    & \qquad
        + \rho(\BA)^{\theta\lambda} \,
            r\bigg(\frac{s-u(\BA)}{\rho(\BA)^\theta}\bigg),
\allowdisplaybreaks
\label{E:DLA}
\\
    & \DC\chi(s|\BA)
\nonumber\\
    & \quad
        = \rho(\BA)^{\theta\lambda}
            \bigg( A_1 \Big( \delta(s-\UA)+\delta(s-\OA) \Big)
                + A_2 \Big( \PV(s-\UA)-\PV(s-\OA) \Big) \bigg)
\nonumber\\
    & \qquad
        + \rho(\BA)^{\theta(\lambda-1)}
            \bigg( A_3 \Big( H(s-\UA)-H(s-\OA) \Big)
                + A_4 \Big( \CI(s-\UA)+\CI(s-\OA) \Big) \bigg)
\nonumber\\
    & \qquad
        + \rho(\BA)^{\theta(\lambda-1)}
            \bigg( -A_42\theta\log\rho(\BA)
                + q\bigg(\frac{s-u(\BA)}{\rho(\BA)^\theta}\bigg)
                    \bigg)
\label{E:DLP}
\end{align}
in the distribution sense in $s$ for all $\BA\in\CC$. Using
\eqref{E:CHISIG} and the product rule \eqref{E:RFV2} we obtain
similar identities for the entropy-flux kernel $\sigma$. For
$\gamma=5/3$ we have $A_2=A_4=0$, so \eqref{E:DLA} and
\eqref{E:DLP} do not contain $\PV$ and $\CI$.

\subsection{Proof of the reduction result}\label{SS:RED}

We essentially follow the arguments in \cite{CL,LPS}. But since we
no longer assume that $\SPT\nu$ is a bounded set, we must ensure
that all terms are indeed well-defined. Let us first fix some
notation.

We choose nonnegative test functions $\varphi, \varphi' \in \D(\R)$
with support in the interval $[-1,1]$ and with integral equal to
one. For $\EPS>0$ we put
$$
    \varphi_\EPS(s) := \EPS^{-1} \varphi(s/\EPS),
    \qquad
    \varphi'_\EPS(s) := \EPS^{-1} \varphi'(s/\EPS)
$$
for all $(s,\EPS)\in\R\times(0,1)$. We then mollify the entropy
kernels: Let
$$
    \chi_\EPS(s|\BA) := \chi(\cdot|\BA)\star\varphi_\EPS(s),
    \qquad\
    \sigma_\EPS(s|\BA) := \sigma(\cdot|\BA)\star\varphi_\EPS(s)
$$
for all $(s,\BA)\in\R\times\CC$, and define $(\chi'_\EPS,
\sigma'_\EPS)$ analogously, using the mollifier $\varphi'_\EPS$
instead. We assume that $\varphi$ and $\varphi'$ are chosen in such
a way that
\begin{equation}
    Z:=\iint_{\R\times\R} H(t-s) \Big(
            \varphi(t)\varphi'(s)-\varphi(s)\varphi'(t)
        \Big) \,ds \,dt
\label{E:NONZ}
\end{equation}
is a positive number. As shown in \cite{CL}, this is always
possible.

The proof of Theorem~\ref{T:CONC} relies on the following two
propositions.

\begin{proposition}\label{P:EINS} There exist a constant
$B>0$ depending on $\lambda$ and the number $Z$ defined in
\eqref{E:NONZ} such that for any nonnegative $\zeta\in\D(\R)$ we
have
\begin{align*}
    & \lim_{\EPS\rightarrow 0}
        \int_\R \Big\LA
                 \DC\chi_\EPS(t) \DC\sigma'_\EPS(t)
                -\DC\sigma_\EPS(t) \DC\chi'_\EPS(t) \Big\RA
            \big\LA\chi(t)\big\RA \zeta(t) \,dt
\\
    &\quad
        = B \int_\CC \rho(\BA)^{1-\theta} \Big(
            \big\LA\chi(\OA)\big\RA \zeta(\OA)
            + \big\LA\chi(\UA)\big\RA \zeta(\UA) \Big)
                \,\nu(d\BA).
\end{align*}
\end{proposition}

\begin{proposition}\label{P:ZWEI} For any test function
$\zeta\in\D(\R)$ we have
\begin{align*}
    & \lim_{\EPS\rightarrow 0} \int_\R \Big\LA
                 \chi(t) \DC\sigma'_\EPS(t)
                -\sigma(t) \DC\chi'_\EPS(t) \Big\RA
            \big\LA\DC\chi_\EPS(t)\big\RA \zeta(t) \,dt
\\
    & \quad = \lim_{\EPS\rightarrow 0} \int_\R \Big\LA
                 \chi(t) \DC\sigma_\EPS(t)
                -\sigma(t) \DC\chi_\EPS(t) \Big\RA
            \big\LA\DC\chi'_\EPS(t) \big\RA \zeta(t) \,dt.
\end{align*}
\end{proposition}

Propositions~\ref{P:EINS} will be proved in Subsection~\ref{SS:E},
Proposition~\ref{P:ZWEI} in Subsection~\ref{SS:Z}. Let us first
show
how they imply Theorem~\ref{T:CONC}. Following the strategy
introduced in \cite{CL} we multiply \eqref{E:COMR} by $\LA
\chi(t)\RA$ and obtain the identity
\begin{align*}
    & \Big\LA\chi(s)\sigma(s')-\sigma(s)\chi(s')\Big\RA
            \big\LA\chi(t)\big\RA
\\
    &\quad = \Big( \big\LA\chi(s)\big\RA
                \big\LA\sigma(s')\big\RA
               -\big\LA\sigma(s)\big\RA
                \big\LA\chi(s')\big\RA \Big)
            \big\LA\chi(t)\big\RA
\end{align*}
for almost all $(s,s',t)\in\R^3$. Cyclic permutation of the
variables yields
\begin{align*}
    & \Big\LA\chi(s')\sigma(t)-\sigma(s')\chi(t)\Big\RA
            \big\LA\chi(s)\big\RA
\\
    &\quad = \Big( \big\LA\chi(s')\big\RA
                 \big\LA\sigma(t)\big\RA
                -\big\LA\sigma(s')\big\RA
                 \big\LA\chi(t)\big\RA \Big)
            \big\LA\chi(s)\big\RA,
\\
    & \Big\LA\chi(t)\sigma(s)-\sigma(t)\chi(s)\Big\RA
            \big\LA\chi(s')\big\RA
\\
    &\quad = \Big( \big\LA\chi(t)\big\RA
                 \big\LA\sigma(s)\big\RA
                -\big\LA\sigma(t)\big\RA
                 \big\LA\chi(s)\big\RA \Big)
            \big\LA\chi(s')\big\RA.
\end{align*}
Summing up all terms, the right-hand sides cancel out, and we find
\begin{align*}
    & \Big\LA\chi(s)\sigma(s')-\sigma(s)\chi(s')\Big\RA
            \big\LA\chi(t)\big\RA
\\
    &\quad
        = \Big\LA\chi(t)\sigma(s')-\sigma(t)\chi(s')\Big\RA
            \big\LA\chi(s)\big\RA
         -\Big\LA\chi(t)\sigma(s)-\sigma(t)\chi(s)\Big\RA
            \big\LA\chi(s')\big\RA.
\end{align*}
We apply the fractional differentiation operator $\DC$ with respect
to $s$ and $s'$, then integrate against the mollifiers
$\varphi_\EPS(t-s)$ and $\varphi'_\EPS(t-s')$ as defined in the
beginning of Subsection~\ref{SS:RED}. Finally, we multiply the
resulting terms by some nonnegative test function $\zeta \in\D(\R)$
and integrate in $t$ over $\R$. Then
\begin{align*}
    & \int_\R \Big\LA
                 \DC\chi_\EPS(t) \DC\sigma'_\EPS(t)
                -\DC\sigma_\EPS(t) \DC\chi'_\EPS(t)
            \Big\RA
            \big\LA\chi(t)\big\RA \zeta(t) \,dt
\\
    &\quad
        = \int_\R \Big\LA
                 \chi(t) \DC\sigma'_\EPS(t)
                -\sigma(t) \DC\chi'_\EPS(t) \Big\RA
            \big\LA\DC\chi_\EPS(t)\big\RA \zeta(t) \,dt
\\
    &\qquad
        -\int_\R \Big\LA
                 \chi(t) \DC\sigma_\EPS(t)
                -\sigma(t) \DC\chi_\EPS(t) \Big\RA
            \big\LA\DC\chi'_\EPS(t) \big\RA \zeta(t) \,dt.
\end{align*}
According to Proposition~\ref{P:EINS}, the right-hand side
converges
to zero as $\EPS\rightarrow 0$ since the two terms have the same
limit. Proposition~\ref{P:ZWEI} describes the limit of the
left-hand
side. Sending $\EPS\rightarrow 0$, we arrive at the identity
\begin{equation}
    B\int_\CC \rho(\BA)^{1-\theta} \Big(
            \big\LA\chi(\OA)\big\RA \zeta(\OA)
            + \big\LA\chi(\UA)\big\RA \zeta(\UA)
                \Big) \,\nu(d\BA) = 0.
\label{E:MNH}
\end{equation}
All terms of the integrand in \eqref{E:MNH} are nonnegative.
Choosing a monotone sequence of $\zeta_k\in\D(\R)$ with
$0\LS\zeta_k\LS 1$ and $\zeta_k \longrightarrow 1$ as
$k\rightarrow\infty$, we get
\begin{equation}
    \int_\CC \rho(\BA)^{1-\theta}
        \LA\chi(\OA)\RA \,\nu(d\BA) = 0,
    \quad\quad
    \int_\CC \rho(\BA)^{1-\theta}
        \LA\chi(\UA)\RA \,\nu(d\BA) = 0,
\label{E:DAS}
\end{equation}
by monotone convergence. Recall that the constant $B$ is strictly
positive.

Consider now the interval $\SC = (\UZ,\OZ)$ defined in
Lemma~\ref{L:CHICON}. If $\SC=\varnothing$, then the representation
\eqref{E:SP} implies that $\SPT\nu\subset V$. If
$\SC\neq\varnothing$, then we find
$$
    \SPT\nu\cap\big\{\BA\in H\colon
        \text{$\OA>\OZ$ or $\UA<\UZ$}\big\} = \varnothing.
$$
%
%
%
Since $\LA\chi(s)\RA>0$ for all $s\in\SC$, from \eqref{E:DAS} and
\eqref{E:SP} we get
$$
    \SPT\nu\cap\big\{\BA\in H\colon \UZ<\UA<\OZ\big\}
        = \varnothing
    \quad\text{and}\quad
    \SPT\nu\cap\big\{\BA\in H\colon \UZ<\OA<\OZ\big\}
        = \varnothing.
$$
%
Therefore the measure $\nu$ must be
contained in the vacuum $V$ and in the isolated point $\BZ :=
(\UZ, \OZ)\in H$. We make an ansatz
$$
    \nu = (1-\omega)\nu_V + \omega\delta_\BZ
    \quad\text{for some $\omega\in[0,1]$,}
$$
where $\nu_V$ is a probability measure supported in the vacuum $V$.
Using this measure in the commutator relation \eqref{E:COMR}, we
find the identity
$$
    (\omega-\omega^2) \Big( -\chi(s|\BZ)\sigma(s'|\BZ)
        +\sigma(s|\BZ)\chi(s'|\BZ) \Big) = 0,
    \qquad\text{a.e. $(s,s')\in\R^2$.}
$$
For some $s,s'\in\SC$ with $s\neq s'$ the second factor does not
vanish, which implies that $\omega\in\{0,1\}$. If $\omega=0$, then
$\nu$ is supported in the vacuum $V$. If $\omega=1$, then $\nu$ is
a Dirac measure at the point $\BZ$. This proves
Theorem~\ref{T:CONC}.


\subsection{Proof of Proposition~\ref{P:EINS}}
\label{SS:E}

As shown in Proposition~\ref{P:EXP}, the fractional differentiation
operator $\DC$ applied to the entropy/entropy flux-kernels creates
distributions such as Dirac measures, principal values, and their
primitives. Up to mollification, the quantities in
Propositions~\ref{P:EINS} and \ref{P:ZWEI} contain products of
these
distributions, so we must carefully argue that all terms are
well-defined.

Let $\varphi_\EPS$, $\varphi'_\EPS$ be the mollifiers from the
beginning of Subsection~\ref{SS:RED} and define
\begin{equation}
    \Phi_\EPS(s,s')
        := \int_\R g(t) \varphi_\EPS(t-s)
            \varphi'_\EPS(t-s') \,dt,
    \quad (s,s')\in\R^2,
\label{E:PPPHI}
\end{equation}
for all $\EPS>0$. Here $g\in C^\alpha(\R)$ is some nonnegative
function with compact support, with $\alpha\in[0,\lambda]$. Now fix
$L>0$ such that $\SPT g \subset B_L(0)$ and define
$$
    B_1:= B_{L+1}(0)
    \quad\text{and}\quad
    B:=B_{L+2}(0).
$$

The proof of Proposition~\ref{P:EINS} is based on the following two
lemmas.

\begin{lemma}
\label{L:EINS}
Let $R$ be a bounded, H\"{o}lder continuous function.
Consider any pair of distributions $T,T'\in\D'(\R)$ from the
following table:
$$
\begin{tabular}{lllll}
    $(T,T')=(\delta,Q)$, &&
    $(T,T')=(\PV,Q)$, &&
    $(T,T')=(Q,Q')$,
\end{tabular}
$$
where $Q,Q'\in\{H,\CI,R\}$. Then there exists a constant $C>0$ such
that
\begin{align}
    & \sup_{\EPS\in(0,1)} \bigg| \iint_{\R\times\R} \Phi_\EPS(s,s')
            \Big[ T(s)T'(s')-T'(s)T(s') \Big] \,ds \,ds' \bigg|
\nonumber\\
    & \quad
        \LS \|g\|_{C^\alpha(\R)}
            \bigg( C \big(1+\|R\|_{C^\alpha(B)}\big)^2 \bigg).
\label{E:BOUND1}
\end{align}
Moreover, we have the following limits:\\
\hspace*{1em} {\em(1)} For $(T,T')=(\delta,H)$ or $(\PV,\CI)$ we
have
\begin{align*}
    & \lim_{\EPS\rightarrow 0}
        \iint_{\R\times\R} \Phi_\EPS(s,s')
            \Big[ \delta(s) H(s') - H(s) \delta(s') \Big]
                \,ds \,ds'
        = Z \,g(0),
\\
    & \lim_{\EPS\rightarrow 0}
        \iint_{\R\times\R} \Phi_\EPS(s,s')
            \Big[ \PV(s) \CI(s') - \CI(s) \PV(s') \Big]
                \,ds \,ds'
        = Z \pi^2 \,g(0).
\intertext{\hspace*{1em} {\em(2)} For all other combinations of $T$
and $T'$ we have}
    & \lim_{\EPS\rightarrow 0}
        \iint_{\R\times\R} \Phi_\EPS(s,s')
            \Big[ T(s) T'(s') - T'(s) T(s') \Big]
                \,ds \,ds'
        = 0.
\end{align*}
The constant $Z>0$ is defined by \eqref{E:NONZ}.
\end{lemma}

{\em Proof.}
Note first that the assumptions on $g$ and on the mollifiers
$\varphi_\EPS$ and $\varphi'_\EPS$ imply that the function
$\Phi_\EPS$
is in $\D(\R\times\R)$. Therefore the pairing
\begin{equation}
    \iint_{\R\times\R} \Phi_\EPS(s,s')
        \Big[ T(s)T'(s')-T'(s)T(s') \Big] \,ds \,ds'
\label{E:TERM1}
\end{equation}
is well-defined for all pairs $(T,T')$ considered. As a function of
$\EPS\in(0,1)$, the integral \eqref{E:TERM1} is smooth. To establish
\eqref{E:BOUND1} it is sufficient to control the behavior as
$\EPS\rightarrow 0$, in which case the singularities become
important.

Note that a substitution of variables yields the identity
\begin{align*}
    & \iint_{\R\times\R} \Phi_\EPS(s,s')
        \Big[ T(s)T'(s')-T'(s)T(s') \Big] \,ds \,ds'
\\
    & \quad
        = \iint_{\R\times\R} M_\EPS(u,u')\,
            \varphi(u) \varphi'(u') \,du \,du',
\end{align*}
where the function $M_\EPS$ is defined as
$$
    M_\EPS(u,u')
        := \int_\R g(t)
            \Big[ T(t-\EPS u) T'(t-\EPS u') - T(t-\EPS u)
                T'(t-\EPS u') \Big] \,dt
$$
for $(u,u')\in\R\times\R$. In the following, we will use the
decomposition \eqref{E:COSINE} of the Cosine Integral into a
logarithm and a H\"{o}lder continuous remainder.
\medskip

{\bf Step~1.} Let $(T,T')=(\delta,H)$. Note that
$$
    \int_\R g(t) \delta(t-\EPS u) H(t-\EPS u') \,dt
        = g(\EPS u) H\big(\EPS(u-u')\big),
$$
with a similar identity if $u$ and $u'$ are interchanged. Therefore
$$
    M_\EPS(u,u')
        = g(\EPS u) H\big(\EPS(u-u')\big)
            - g(\EPS u')  H\big(\EPS(u'-u)\big),
$$
which implies the estimate
$$
    |M_\EPS(u,u')|
        \LS 2\|g\|_{L^\infty(\R)}.
$$
By dominated convergence, we obtain
\begin{align*}
    & \lim_{\EPS\rightarrow 0} \iint_{\R\times\R} \Phi_\EPS(s,s')
        \Big[ \delta(s)H(s')-H(s)\delta(s') \Big] \,ds \,ds'
\\
    & \quad
        = g(0) \bigg( \iint_{\R\times\R} \Big[ H(u-u')-H(u'-u)
\Big]
            \varphi(u) \varphi'(u') \,du \,du' \bigg).
\end{align*}
The integral on the right-hand side coincides with $Z>0$ defined in
\eqref{E:NONZ}.
\medskip

{\bf Step~2.} Let $(T,T')=(\delta,\log|\cdot|)$. Note that
$$
    \int_\R g(t) \delta(t-\EPS u) \log|t-\EPS u'| \,dt
        = g(\EPS u) \log|\EPS(u-u')|,
$$
with a similar identity if $u$ and $u'$ are interchanged. Therefore
$$
    M_\EPS(u,u')
        = \Big[ g(\EPS u)-g(\EPS u') \Big] \log|\EPS(u-u')|.
$$
We obtain the estimate
\begin{equation}
    |M_\EPS(u,u')|
        \LS \|g\|_{C^\alpha(\R)} \Big( \big(\EPS|u-u'|\big)^\alpha
            \big|\log|\EPS(u-u')|\big| \Big).
\label{E:QWEINS}
\end{equation}
Since the supports of $\varphi$ and $\varphi'$ are contained in
$[-1,1]$, the right-hand side of \eqref{E:QWEINS} is
uniformly bounded and converges to zero as $\EPS\rightarrow 0$,
yielding
$$
    \lim_{\EPS\rightarrow 0} \iint_{\R\times\R} \Phi_\EPS(s,s')
        \Big[ \delta(s)\log|s'|-\log|s|\delta(s') \Big] \,ds \,ds'
            = 0.
$$
\smallskip

{\bf Step~3.} Let $(T,T')=(\delta,R)$. Note that
$$
    \int_\R g(t) \delta(t-\EPS u) R(t-\EPS u') \,dt
        = g(\EPS u) R\big(\EPS(u-u')\big),
$$
with a similar identity if $u$ and $u'$ are interchanged. Therefore
$$
    M_\EPS(u,u')
        = g(\EPS u) R\big(\EPS(u-u')\big)
            - g(\EPS u') R\big(\EPS(u'-u)\big),
$$
which implies the estimate
$$
    |M_\EPS(u,u')|
        \LS \|g\|_{L^\infty(\R)}
            \Big( 2\|R\|_{L^\infty(\R)} \Big).
$$
By dominated convergence, we then obtain
$$
    \lim_{\EPS\rightarrow 0} \iint_{\R\times\R} \Phi_\EPS(s,s')
        \Big[ \delta(s)R(s')-R(s)\delta(s') \Big] \,ds \,ds'
            = 0.
$$
\smallskip

{\bf Step~4.} Let $(T,T')=(\PV,H)$. A substitution of variables
yields
$$
    \int_\R g(t) \PV(t-\EPS u) H(t-\EPS u') \,dt
        = \int_{-\EPS(u-u')}^\infty \PV(s) g(s+\EPS u) \,ds,
$$
with a similar identity if $u$ and $u'$ are interchanged. Therefore
$$
    M_\EPS(u,u')
        = \int_{-\EPS(u-u')}^\infty \PV(s) g(s+\EPS u) \,ds
            -\int_{-\EPS(u'-u)}^\infty \PV(s) g(s+\EPS u') \,ds.
$$
Let us assume that $u>u'$, the converse case being similar. We
decompose
\begin{align*}
    & \int_{-\EPS(u-u')}^\infty \PV(s) g(s+\EPS u) \,ds
\\
    & \quad
        = \int_{-\EPS(u-u')}^{\EPS(u-u')} \PV(s) g(s+\EPS u) \,ds
            + \int_{\EPS(u-u')}^\infty \PV(s) g(s+\EPS u) \,ds.
\end{align*}
By symmetry, the first integral on the right-hand side can be
rewritten as
$$
    \int_{-\EPS(u-u')}^{\EPS(u-u')} \PV(s) g(s+\EPS u) \,ds
        = \int_{-\EPS(u-u')}^{\EPS(u-u')} \PV(s) \Big[ g(s+\EPS u)
            -g(\EPS u) \Big ] \,ds,
$$
which implies the estimate
\begin{align*}
    \bigg| \int_{-\EPS(u-u')}^{\EPS(u-u')} \PV(s) g(s+\EPS u)
            \,ds \bigg|
        & \LS \|g\|_{C^\alpha(\R)} \int_{-\EPS(u-u')}^{\EPS(u-u')}
            |s|^{\alpha-1} \,ds
\\[1ex]
        & = \|g\|_{C^\alpha(\R)} \Big( 2 \alpha^{-1}
            \big(\EPS|u-u'|\big)^\alpha \Big).
\end{align*}
The right-hand side is uniformly bounded and vanishes as
$\EPS\rightarrow 0$. Now
\begin{align*}
    & \int_{\EPS(u-u')}^\infty \PV(s) g(s+\EPS u) \,ds
        -\int_{-\EPS(u'-u)}^\infty \PV(s) g(s+\EPS u') \,ds
\\
    & \quad
        = \int_{\EPS(u-u')}^\infty \PV(s) \Big[ g(s+\EPS u)
            - g(s+\EPS u') \Big] \,ds,
\end{align*}
which implies the estimate
\begin{align}
    & \bigg| \int_{\EPS(u-u')}^\infty \PV(s) g(s+\EPS u) \,ds
        -\int_{-\EPS(u'-u)}^\infty \PV(s) g(s+\EPS u') \,ds
            \bigg|
\nonumber\\
    & \quad
        \LS \|g\|_{C^\alpha(\R)} \Bigg( \big(\EPS|u-u'|
                \big)^\alpha
            \int_{\EPS(u-u')}^{L+1} \frac{ds}{s} \Bigg)
\nonumber\\
    & \quad
        = \|g\|_{C^\alpha(\R)} \bigg(\big (\EPS|u-u'|
                \big)^\alpha
            \Big[ \log(L+1)-\log|\EPS(u-u')| \Big] \bigg).
\label{E:QWZWEI}
\end{align}
Recall that $\SPT g\subset B_L(0)$. The right-hand side of
\eqref{E:QWZWEI} is uniformly bounded and converges to zero as
$\EPS\rightarrow 0$. Combining the above estimates we get
$$
    \lim_{\EPS\rightarrow 0} \iint_{\R\times\R} \Phi_\EPS(s,s')
        \Big[ \PV(s)H(s')-H(s)\PV(s') \Big] \,ds \,ds'
            = 0.
$$
\smallskip

{\bf Step~5.} Let $(T,T')=(\PV,\log|\cdot|)$. A substitution of
variables yields
$$
    \int_\R g(t) \PV(t-\EPS u) \log|t-\EPS u'| \,dt
        = \int_\R \PV(s) g(s+\EPS u) \log|s+\EPS(u-u')| \,ds,
$$
with a similar identity if $u$ and $u'$ are interchanged. We now
decompose
\begin{align}
    & M_\EPS(u,u')
\nonumber\\
    & \quad
        = \int_{B_1} \PV(s) \Big[
            g(\EPS u') \log|s+\EPS(u-u')|
                -g(\EPS u) \log|s+\EPS(u'-u)| \Big] \,ds
\nonumber\\
    & \qquad
        +\int_{B_1} \PV(s) \Big[
            \big( g(s+\EPS u)-g(\EPS u')
                \big) \log|s+\EPS(u-u')| \Big] \,ds
\nonumber\\
    & \qquad
        -\int_{B_1} \PV(s) \Big[
            \big( g(s+\EPS u')-g(\EPS u)
                \big) \log|s+\EPS(u'-u)| \Big] \,ds.
\label{E:QWDREI}
\end{align}
Note that the function
$$
    \zeta_a(t) := \big( g(t+a)-g(a) \big) \log|t|,
    \quad t\in\R,
$$
is H\"{o}lder continuous for all $a\in\R$. Therefore we can
estimate
\begin{align*}
    & \bigg| \int_{B_1} \PV(s) \Big[ \big( g(s+\EPS u)-g(\EPS u')
            \big) \log|s+\EPS(u-u')| \Big] \,ds \bigg|
\\
    & \quad
        = \bigg| \int_{B_1} \PV(s) \Big[ \zeta_{\EPS u'}\big(s+\EPS
            (u-u')\big)-\zeta_{\EPS u'}\big(\EPS(u-u')\big)
                \Big] \,ds \bigg|
\\
    & \quad
        \LS \|g\|_{C^\alpha(\R)} \bigg( C\int_B |s|^{\alpha'-1}
            \,ds \bigg),
\end{align*}
with $\alpha'<\alpha$ and $C>0$ some constant. Moreover,
we find
\begin{align*}
    & \lim_{\EPS\rightarrow 0} \int_{B_1} \PV(s) \Big[
        \big( g(s+\EPS u)-g(\EPS u') \big)
            \log|s+\EPS(u-u')| \Big] \,ds
\\
    & \quad
        = \int_{B_1} \PV(s) \Big[ g(s)-g(0) \Big] \log|s| \,ds.
\end{align*}
The same reasoning applies with $u$ and $u'$ interchanged, with the
same limit. Therefore the last two terms in \eqref{E:QWDREI} are
bounded and vanish as $\EPS\rightarrow 0$.

To control the first term on the right-hand side of
\eqref{E:QWDREI}, we write
\begin{align*}
    \int_{B_1} \PV(s) \log|s+\EPS(u-u')| \,ds
        & = \int_{-\EPS(u-u')}^{\EPS(u-u')} \PV(s)
            \log|s+\EPS(u-u')| \,ds
\nonumber\\
        & \quad
            + \int_{\EPS(u-u')}^{L+1} \PV(s) \log\bigg|
                \frac{s+\EPS(u-u')}{s-\EPS(u-u')}\bigg| \,ds,
\end{align*}
assuming without loss of generality that $u-u'>0$. Now we have
\begin{align*}
    \int_{-\EPS(u-u')}^{\EPS(u-u')} \PV(s)
            \log|s+\EPS(u-u')| \,ds
        & = \pi^2/4,
\\
    \int_{\EPS(u-u')}^{L+1} \PV(s) \log\bigg|
            \frac{s+\EPS(u-u')}{s-\EPS(u-u')}\bigg| \,ds
        & = \pi^2/4 - h\big(\EPS(u-u')\big),
\end{align*}
where $h$ is a smooth, increasing function with $\lim_{s\rightarrow
0}h(s)=0$. If $u$ and $u'$ are interchanged, we obtain the same
quantities with a minus sign. Therefore
\begin{align*}
    & \int_{B_1} \PV(s) \Big[
        g(\EPS u') \log|s+\EPS(u-u')|
            -g(\EPS u) \log|s+\EPS(u'-u)| \Big] \,ds
\\
    & \quad
        = \big( g(\EPS u)+g(\EPS u') \big)
            \Big( \pi^2/2-h\big(\EPS(u-u')\big) \Big).
\end{align*}
This left-hand side is bounded in absolute value by $\pi^2
\|g\|_{L^\infty (\R)}$ and converges to the limit $\pi^2 g(0)$.
Combining all estimates, we conclude that
$$
    |M_\EPS(u,u')| \LS C \|g\|_{C^\alpha(\R)},
$$
with $C>0$ some constant. By dominated convergence, we find
\begin{align*}
    & \lim_{\EPS\rightarrow 0} \iint_{\R\times\R} \Phi_\EPS(s,s')
        \Big[ \PV(s)\log|s'|-\log|s|\PV(s') \Big] \,ds \,ds'
\\
    & \quad
        = g(0) \bigg( \pi^2 \iint_{\R\times\R}
            \Big[ H(u-u')-H(u'-u) \Big]
                \varphi(u) \varphi'(u') \,du \,du' \bigg).
\end{align*}
The integral on the right-hand side coincides with $Z>0$
defined in \eqref{E:NONZ}.
\medskip

{\bf Step~6.} Let $(T,T')=(\PV,R)$. A substitution of variables
yields
$$
    \int_\R g(t) \PV(t-\EPS u) R(t-\EPS u') \,dt
        = \int_\R \PV(s) g(s+\EPS u) R\big(s+\EPS(u-u')\big)
            \,ds,
$$
with a similar identity if $u$ and $u'$ are interchanged. Therefore
\begin{align*}
    & M_\EPS(u,u')
\\
    & \quad
        = \int_\R \PV(s) \Big[
            g(s+\EPS u)R\big(s+\EPS(u-u')\big)
                - g(s+\EPS u')R\big(s+\EPS(u'-u)\big) \Big] \,ds.
\end{align*}
Since $g$ and $R$ are H\"{o}lder continuous functions, we can
estimate
\begin{align*}
    & \bigg| \int_\R \PV(s) g(s+\EPS u) R\big(s+\EPS(u-u')\big)
        \,ds \bigg|
\\
    & \quad
        = \bigg| \int_{B_1} \PV(s) \Big[
            g(s+\EPS u) R\big(s+\EPS(u-u')\big)
            -g(\EPS u) R\big(\EPS(u-u')\big) \Big] \,ds \bigg|
\\
    & \quad
        \LS \|g\|_{C^\alpha(\R)} \bigg( \|R\|_{C^\alpha(\R)}
            \int_B |s|^{\alpha-1} \,ds \bigg).
\end{align*}
By dominated convergence, we then have
$$
    \lim_{\EPS\rightarrow 0}
            \int_\R g(t) \PV(t-\EPS u) R(t-\EPS u')\,dt
        = \int_B \PV(s) \Big[ g(s)R(s)-g(0)R(0) \Big] \,ds.
$$
The same reasoning applies with $u$ and $u'$ interchanged. We
obtain
the estimate
$$
    |M_\EPS(u,u')|
        \LS \|g\|_{C^\alpha(\R)} \Big( C\|R\|_{C^\alpha(\R)} \Big)
$$
with $C>0$ some constant, and the convergence
$$
    \lim_{\EPS\rightarrow 0} \iint_{\R\times\R} \Phi_\EPS(s,s')
        \Big[ \PV(s)R(s')-R(s)\PV(s') \Big] \,ds \,ds'
            = 0.
$$
\medskip

{\bf Step~7.} Finally, let $(T,T')=(Q,Q')$ with $Q,Q'\in\{
H,\log|\cdot|,R\}$. We have
\begin{align*}
    |M_\EPS(u,u')|
    & \LS \|g\|_{L^\infty(\R)} \Big(
        \|Q(\cdot-\EPS u)-Q(\cdot-\EPS u')\|_{L^2(B)}
            \|Q'(\cdot-\EPS u')\|_{L^2(B)}
\\
    & \qquad
        +\|Q(\cdot-\EPS u')\|_{L^2(B)}
            \|Q'(\cdot-\EPS u')-Q'(\cdot-\EPS u)\|_{L^2(B)}
        \Big).
\end{align*}
Since $Q,Q'\in W^{\beta,2}_\LOC(\R)$ for all $\beta<1$, the
right-hand side is uniformly bounded and converges to zero as $\EPS
\rightarrow 0$. By dominated convergence, we get that
$$
    \lim_{\EPS\rightarrow 0} \iint_{\R\times\R} \Phi_\EPS(s,s')
        \Big[ Q(s)Q'(s')-Q'(s)Q(s') \Big] \,ds \,ds'
            = 0.
$$
The proof of Lemma~\ref{P:EINS} is now complete.
\qed

\begin{lemma}
\label{L:ZWEI}
Let $R$ be a bounded, H\"{o}lder continuous function.
Consider any pair of distributions $T,T'\in\D'(\R)$ from the
following table:
$$
\begin{tabular}{lllll}
    $\{T,T'\}=\{\delta,\delta\}$, &&
    $\{T,T'\}=\{\PV,\PV\}$, &&
    $\{T,T'\}=\{Q,Q\}$,
\\
    $\{T,T'\}=\{\delta,\PV\}$, &&
    $\{T,T'\}=\{\PV,Q\}$,
\\
    $\{T,T'\}=\{\delta,Q\}$, &&
\end{tabular}
$$
where $Q\in\{H,\CI,R\}$. Then there exists a constant $C>0$ such
that
\begin{align}
    & \sup_{\EPS\in(0,1)} \bigg| \iint_{\R\times\R}
(s-s')\Phi_\EPS(s,s')
        \Big[ T(s) T'(s') \Big] \,ds \,ds' \bigg|
\nonumber\\
    & \quad
        \LS \|g\|_{C^\alpha(\R)}
            \bigg( C \big(1+\|R\|_{C^\alpha(B)}\big)^2 \bigg).
\label{E:BOUND2}
\end{align}
Moreover, we have the following limits:\\
\hspace*{1em} {\em(1)} For $\{T,T'\}=\{\delta,\PV\}$ we have
\begin{align}
    & \lim_{\EPS\rightarrow 0}
        \iint_{\R\times\R} (s-s') \Phi_\EPS(s,s')
            \Big[ \PV(s) \delta(s') + \delta(s) \PV(s')
                \Big] \,ds \,ds'
        = 0.
\label{E:EXTIN}
\intertext{\hspace*{1em} {\em(2)} For all other combinations of $T$
and $T'$ we have}
    & \lim_{\EPS\rightarrow 0}
        \iint_{\R\times\R} (s-s') \Phi_\EPS(s,s')
            \Big[ T(s) T'(s') \Big]
                \,ds \,ds'
        = 0.
\nonumber
\end{align}
\end{lemma}

{\em Proof.} Note first that the map $(s,s')\mapsto
(s-s')\Phi_\EPS(s,s')$ is in $\D(\R\times\R)$ since the function
$\Phi_\EPS$ is smooth with compact support. This follows from
\eqref{E:PPPHI}, and from the assumptions on $g$ and $\varphi_\EPS$,
$\varphi'_\EPS$. Therefore the pairing with products of
distributions is well-defined. As in the proof of
Lemma~\ref{L:EINS}, in order to establish the bound \eqref{E:BOUND2}
it is sufficient to consider the behavior as $\EPS\rightarrow 0$.
\medskip

{\bf Step~1.} We immediately find that
$$
    \iint_{\R\times\R} (s-s')\Phi_\EPS(s,s')
            \Big[ \delta(s) \delta(s') \Big] \,ds \,ds'
        = 0.
$$
\smallskip

{\bf Step~2.} We have the identity
\begin{align}
    & \iint_{\R\times\R} (s-s')\Phi_\EPS(s,s')
            \Big[ \PV(s) \delta(s') \Big] \,ds \,ds'
\nonumber\\
    & \quad
        = \int_\R g(t) \bigg( \int_\R \varphi_\EPS(t-s)
             \Big[ s\PV(s) \Big] \,ds \bigg)
                \varphi'_\EPS(t) \,dt
\nonumber\\
    & \quad
        = \int_\R g(t) \varphi'_\EPS(t) \,dt,
\label{E:DEPR}
\end{align}
where we used the fact that $s\PV(s)=1$. We can therefore estimate
$$
    \bigg| \iint_{\R\times\R} (s-s')\Phi_\EPS(s,s')
            \Big[ \PV(s) \delta(s') \Big] \,ds \,ds' \bigg|
        \LS \|g\|_{L^\infty(\R)}.
$$
Moreover, by continuity of $g$ we obtain the convergence
$$
    \lim_{\EPS\rightarrow 0} \iint_{\R\times\R}
        (s-s')\Phi_\EPS(s,s')
            \Big[ \PV(s) \delta(s') \Big] \,ds \,ds'
                = g(0).
$$
If we reverse the order of the distributions, the same reasoning
applies. The resulting term converges to $-g(0)$ as
$\EPS\rightarrow
0$, so the claim \eqref{E:EXTIN} follows.
\medskip

{\bf Step~3.} We have the identity
\begin{align*}
    & \iint_{\R\times\R} (s-s')\Phi_\EPS(s,s')
            \Big[ \log|s| \delta(s') \Big] \,ds \,ds'
\\
    & \quad
        = \int_\R g(t) \bigg( \int_\R
            \varphi_\EPS(t-s) \Big[ s\log|s| \Big] \,ds
                \bigg) \varphi'_\EPS(t) \,dt.
\end{align*}
We can therefore estimate as follows:
$$
    \bigg| \iint_{\R\times\R} (s-s')\Phi_\EPS(s,s')
            \Big[ \log|s| \delta(s') \Big] \,ds \,ds' \bigg|
        \LS \|g\|_{L^\infty(\R)} \bigg(
            \sup_{|s|\LS 2\EPS} \big|s\log|s|\big| \bigg).
$$
The right-hand side converges to zero as $\EPS\rightarrow 0$.
Similar
reasoning applies if the function $\log|\cdot|$ is replaced by $H$
or $R$, and if the order of the distributions are reversed. In
particular, we have the estimate
\begin{equation}
    \bigg| \iint_{\R\times\R} (s-s')\Phi_\EPS(s,s')
            \Big[ R(s) \delta(s') \Big] \,ds \,ds' \bigg|
        \LS \|g\|_{L^\infty(\R)} \bigg(
            2\EPS \|R\|_{L^\infty(B)} \bigg),
\label{E:DEHOE}
\end{equation}
which again vanishes in the limit $\EPS\rightarrow 0$.
\medskip

{\bf Step~4.} We have the identity
\begin{align*}
    & \iint_{\R\times\R} (s-s')\Phi_\EPS(s,s')
            \Big[ \PV(s) \PV(s') \Big] \,ds \,ds'
\\
    & \quad
        = \int_\R g(t) \Bigg\{
            \bigg( \int_\R \varphi_\EPS(t-s)
                \Big[ s\PV(s) \Big] \,ds \bigg)
            \bigg( \int_\R \varphi'_\EPS(t-s')
                \PV(s') \,ds' \bigg)
\\
    & \quad\qquad
        - \bigg( \int_\R \varphi_\EPS(t-s)
                \PV(s) \,ds \bigg)
            \bigg( \int_\R \varphi'_\EPS(t-s')
                \Big[ s'\PV(s') \Big] \,ds' \bigg)
        \Bigg\} \,dt
\\
    & \quad
        = \iint_{\R\times \R} g(t)
            \Big[ \varphi'_\EPS(t-s)-\varphi_\EPS(t-s)
                \Big] \PV(s) \,ds \,dt,
\end{align*}
where we used that $s\PV(s)=1$. After a substitution of
variables, we get
\begin{align}
    & \iint_{\R\times\R} (s-s')\Phi_\EPS(s,s')
            \Big[ \PV(s) \PV(s') \Big] \,ds \,ds'
\nonumber\\
    & \quad
        = \iint_{\R\times \R} g(s+w)
            \Big[ \varphi'_\EPS(w)-\varphi_\EPS(w)
                \Big] \PV(s) \,ds \,dw
\nonumber\\
    & \quad
        = \int_\R \Big[ \varphi'_\EPS(w)-\varphi_\EPS(w)
            \Big] \bigg(
                \int_{B_1} \Big[ g(s+w)-g(w) \Big]
                    \PV(s) \,ds \bigg) \,dw.
\label{E:OLK}
\end{align}
Now we estimate
\begin{align}
    & \bigg| \iint_{\R\times\R} (s-s')\Phi_\EPS(s,s')
            \Big[ \PV(s) \PV(s') \Big] \,ds \,ds' \bigg|
\nonumber\\
    & \quad
        \LS \|g\|_{C^\alpha(\R)} \bigg(
            2 \int_B |t|^{\alpha-1} \,dt \bigg).
\label{E:SUPE}
\end{align}
Note that the map
$$
    \zeta(w) :=
        \int_{B_1} \Big[ g(s+w)-g(w) \Big]
            \PV(w) \,dw
$$
is H\"{o}lder continuous and locally bounded. Therefore we obtain
\begin{align*}
    & \lim_{\EPS\rightarrow 0} \int_\R \varphi_\EPS(w) \bigg(
        \int_{B_1} \Big[ g(s+w)-g(w) \Big]
                    \PV(s) \,ds \bigg) \,dw
\\
    & \quad
        = \int_B \Big[ g(s)-g(0) \Big] \PV(s) \,ds.
\end{align*}
The same holds with $\varphi'_\EPS$ in place of
$\varphi_\EPS$, therefore \eqref{E:OLK} converges to zero.
\medskip

{\bf Step~5.} We have the identity
\begin{align}
    & \iint_{\R\times\R} (s-s')\Phi_\EPS(s,s')
            \Big[ \log|s| \PV(s') \Big] \,ds \,ds'
\nonumber\\
    & \quad
        = \int_\R g(t)
            \bigg( \int_\R \varphi_\EPS(t-s)
                \Big[ s\log|s| \Big] \,ds \bigg)
            \bigg( \int_\R \varphi'_\EPS(t-s')
                \PV(s') \,ds' \bigg) \,dt
\nonumber\\
    & \qquad\quad
        - \int_\R g(t)
            \bigg( \int_\R \varphi_\EPS(t-s)
                \log|s| \,ds \bigg) \,dt,
\label{E:POI}
\end{align}
where we used that $s'\PV(s')=1$. The second term can be
estimated as
\begin{align}
    \bigg| \int_\R g(t) \bigg( \int_\R \varphi_\EPS(t-s)
            \log|s| \,ds \bigg) \,dt \bigg|
    & = \bigg| \int_\R \varphi_\EPS(w) \bigg( \int_\R g(t)
            \log|t-w| \,dt \bigg) \,dw \bigg|
\nonumber\\
    & \LS \|g\|_{L^\infty(\R)} \bigg(
            \int_B \big|\log|t|\big| \,dt \bigg).
\label{E:LMN}
\end{align}
As in Step~4 we find that the map
$$
    w \mapsto \int_\R g(t) \log|t-w| \,dt
$$
is H\"{o}lder continuous and locally bounded, which implies that
\begin{equation}
    \lim_{\EPS\rightarrow 0} \int_\R g(t)
        \bigg( \int_\R \varphi_\EPS(t-s) \log|s| \,ds \bigg) \,dt
            = \int_\R g(t) \log|t| \,dt.
\label{E:DFG}
\end{equation}
For the first term in \eqref{E:POI} we argue as follows: We
introduce the function
\begin{equation}
    \zeta_\EPS(s') := \int_\R \Bigg( g(t) \int_\R
        \varphi_\EPS(t-s) \Big[ s\log|s| \Big] \,ds \Bigg)
            \varphi'_\EPS(t-s') \,dt
\label{E:HOE}
\end{equation}
for all $s'\in\R$. Since $s\mapsto s\log|s|$ is H\"{o}lder
continuous for all H\"{o}lder exponents less than one, we find that
$\zeta_\EPS$ converges strongly in the $C^\alpha(\R)$-norm to
\begin{equation}
    \zeta(s') := g(s') \Big[s'\log|s'|\Big],
    \quad s'\in\R.
\label{E:LIZE}
\end{equation}
In particular, the $C^\alpha(\R)$-norm of $\zeta_\EPS$ is bounded
uniformly in $\EPS\in (0,1)$, and can in fact be estimated by
$C\|g\|_{C^\alpha(\R)}$, with $C>0$ some constant. Hence
\begin{align}
    & \bigg| \int_\R g(t)
            \bigg( \int_\R \varphi_\EPS(t-s)
                \Big[ s\log|s| \Big] \,ds \bigg)
            \bigg( \int_\R \varphi'_\EPS(t-s')
                \PV(s') \,ds' \bigg) \,dt \bigg|
\nonumber\\
    & \quad
        = \bigg| \int_B \PV(s') \Big[ \zeta_\EPS(s')-\zeta_\EPS(0)
            \Big] \,ds' \bigg|
\nonumber\\
    & \quad
        \LS \|g\|_{C^\alpha(\R)} \bigg( C \int_B
                |s'|^{\alpha-1} \,ds' \bigg).
\label{E:ERF}
\end{align}
From the strong convergence of $\zeta_\EPS$ in the H\"{o}lder-norm
we obtain
\begin{align}
    & \lim_{\EPS\rightarrow 0} \int_\R g(t)
            \bigg( \int_\R \varphi_\EPS(t-s)
                \Big[ s\log|s| \Big] \,ds \bigg)
            \bigg( \int_\R \varphi'_\EPS(t-s')
                \PV(s') \,ds' \bigg) \,dt
\nonumber\\
    & \quad
        = \int_B \PV(s') \Big[\zeta(s')-\zeta(0)\Big] \,ds'
\nonumber\\
    & \quad
        = \int_B g(s') \log|s'| \,ds',
\label{E:CVB}
\end{align}
using that $s'\PV(s')=1$ and $\zeta(0)=0$. Because of \eqref{E:DFG}
and \eqref{E:CVB}, the right-hand side of \eqref{E:POI} vanishes as
$\EPS\rightarrow 0$. The same holds with $\log|\cdot|$
replaced by $H$ or $R$, and with the order of the distributions
reversed. We have
\begin{align*}
    & \bigg| \iint_{\R\times\R} (s-s')\Phi_\EPS(s,s')
            \Big[ R(s) \PV(s') \Big] \,ds \,ds' \bigg|
\\
    & \quad
        \LS \|g\|_{L^\infty(\R)} \bigg(
                \|R\|_{L^1(B)} \bigg)
            + \|g\|_{C^\alpha(\R)} \bigg( C
                \|R\|_{C^\alpha(B)}
                \int_B |s'|^{\alpha-1} \,ds' \bigg),
\end{align*}
which implies the desired estimate.
\medskip

{\bf Step~6.} We have the identity
\begin{align}
    & \iint_{\R\times\R} (s-s')\Phi_\EPS(s,s')
            \Big[ \log|s| \log|s'| \Big] \,ds \,ds'
\nonumber\\
    & \quad
        = \int_\R g(t) \Bigg\{
            \bigg( \int_\R \varphi_\EPS(t-s)
                \Big[ s\log|s| \Big] \,ds \bigg)
            \bigg( \int_\R \varphi'_\EPS(t-s')
                \log|s'| \,ds' \bigg)
\nonumber\\
    & \qquad
        - \bigg( \int_\R \varphi_\EPS(t-s)
                \log|s| \,ds \bigg)
            \bigg( \int_\R \varphi'_\EPS(t-s')
                \Big[ s'\log|s'| \Big] \,ds' \bigg)
        \Bigg\} \,dt.
\label{E:LOGLOG}
\end{align}
Using again the function $\zeta_\EPS$ defined in \eqref{E:HOE},
which converges strongly in the $\sup$-norm to the limit
\eqref{E:LIZE}, we can now estimate
\begin{align*}
    & \bigg| \int_\R g(t)
            \bigg( \int_\R \varphi_\EPS(t-s)
                \Big[ s\log|s| \Big] \,ds \bigg)
            \bigg( \int_\R \varphi'_\EPS(t-s')
                \log|s'| \,ds' \bigg) \,dt \bigg|
\\
    & \quad
        = \bigg| \int_B \log|s'| \zeta_\EPS(s') \,ds' \bigg|
\\
    & \quad
        \LS \|g\|_{L^\infty(\R)} \bigg( C \int_B
                \log|s'| \,ds' \bigg),
\end{align*}
with $C>0$ some constant. From the strong convergence of
$\zeta_\EPS$, we obtain
\begin{align*}
    & \lim_{\EPS\rightarrow 0} \int_\R g(t)
            \bigg( \int_\R \varphi_\EPS(t-s)
                \Big[ s\log|s| \Big] \,ds \bigg)
            \bigg( \int_\R \varphi'_\EPS(t-s')
                \PV(s') \,ds' \bigg) \,dt
\\
    & \quad
        = \int_B \log|s'| \zeta(s') \,ds'
\\
    & \quad
        = \int_B g(s') s'\big(\log|s'|)^2 \,ds'.
\end{align*}
The same limit is obtained with primed and unprimed terms
interchanged, so the left-hand side of \eqref{E:LOGLOG} vanishes as
$\EPS\rightarrow 0$. Any other combination of functions from
$\{\log|\cdot|,H,R\}$ can be handled in the same way. We have
$$
    \bigg| \iint_{\R\times\R} (s-s')\Phi_\EPS(s,s')
            \Big[ R(s) R(s') \Big] \,ds \,ds' \bigg|
        \LS \|g\|_{L^\infty(\R)} \bigg(
            2\|R\|_{L^\infty(B)} \|R\|_{L^1(B)} \bigg),
$$
with similar estimates for the remaining combinations.
\qed

{\em Proof of Proposition~\ref{P:EINS}.} Using \eqref{E:CHISIG} and
\eqref{E:RFV2} we find the identity
\begin{align}
    & \DC\chi(s|\BA)\, \DC\sigma(s'|\BA)
        -\DC\sigma(s|\BA)\, \DC\chi(s'|\BA)
\nonumber\\
    & \vphantom{\Big[}\quad
        = \theta (s'-s)\, \DC\chi(s|\BA)\, \DC\chi(s'|\BA)
\nonumber\\
    & \qquad
        + \theta(\lambda+1) \Big[ \DC\chi(s|\BA)\,
                \DS\chi(s'|\BA)
            -\DS\chi(s|\BA)\, \DC\chi(s'|\BA) \Big],
\label{E:CSMSC2}
\end{align}
which holds distributionally in $(s,s')\in\R\times\R$ for all
$\BA\in\CC$. Let us consider the first term on the right-hand side.
We fix some $\BA\in\CC$ and integrate against the function
\eqref{E:PPPHI}. We then want to use the expansion \eqref{E:DLP} to
control
\begin{equation}
    \iint_{\R\times\R} (s-s')\Phi_\EPS(s,s') \Big[
        \DC\chi(s|\BA)\, \DC\chi(s'|\BA) \Big] \,ds \,ds'.
\label{E:INTEE}
\end{equation}
Note that $\DC\chi(s|\BA)$ is singular at $s=\UA$ and
$s=\OA$, and smooth otherwise. A straightforward, but tedious
application of Proposition~\ref{P:ZWEI} shows
\begin{align}
    & \sup_{\EPS\in(0,1)} \bigg| \iint_{\R\times\R}
        (s-s')\Phi_\EPS(s,s') \Big[ \DC\chi(s|\BA)\,
            \DC\chi(s'|\BA) \Big] \,ds \,ds' \bigg|
\nonumber\\
    & \quad
        \LS \|g\|_{C^\alpha(R)} \bigg\{
            C \rho(\BA)^{2\theta\lambda}
                \Big( 1+\rho(\BA)^{-\theta} \Big)
                \Big( 1 + \rho(\BA)^{-\alpha\theta}
                    + |\log\rho(\BA)| \Big) \bigg\},
\label{E:UNIB}
\end{align}
with $C>0$ some constant independent of $\BA$. Since $2\lambda-1>0$
for $\gamma\in(1,5/3]$, the right-hand side of \eqref{E:UNIB}
vanishes as $\rho(\BA)\rightarrow 0$, if $\alpha$ is chosen small
enough. For $\rho(\BA)$ large, \eqref{E:UNIB} grows at most
linearly because $2\theta\lambda = 1-\theta<1$. By
Proposition~\ref{P:ZWEI} and the dominated convergence theorem, we
obtain
$$
    \lim_{\EPS\rightarrow 0}
        \iint_{\R\times\R} (s-s')\Phi_\EPS(s,s') \Big\LA
            \DC\chi(s)\, \DC\chi(s') \Big\RA \,ds \,ds'
        = 0.
$$
For the second term in \eqref{E:CSMSC2} we argue similarly: Again we
have a bound
\begin{align*}
    & \sup_{\EPS\in(0,1)} \bigg| \iint_{\R\times\R} \Phi_\EPS(s,s')
        \Big[ \DC\chi(s|\BA)\, \DS\chi(s'|\BA)
            -\DS\chi(s|\BA)\, \DC\chi(s'|\BA) \Big] \,ds \,ds'
\bigg|
\nonumber\\
    & \quad
        \LS \|g\|_{C^\alpha(R)} \bigg\{
            C \rho(\BA)^{2\theta\lambda}
                \Big( 1+\rho(\BA)^{-\theta} \Big)
                \Big( 1 + \rho(\BA)^{-\alpha\theta}
                    + \big|\log\rho(\BA)\big| \Big) \bigg\}
\end{align*}
with $C>0$ some constant, as follows from the expansions
\eqref{E:DLA} and \eqref{E:DLP}. We use Proposition~\ref{P:EINS}
and the dominated convergence theorem to obtain
\begin{align*}
    & \lim_{\EPS\rightarrow 0}
        \iint_{\R\times\R} \Phi_\EPS(s,s') \Big\LA
            \DC\chi(s)\, \DS\chi(s')
                -\DS\chi(s)\, \DC\chi(s')
            \Big\RA \,ds \,ds'
\\
    & \quad
        = (A_1^2+\pi^2 A_2^2)Z \int_\CC
            \rho(\BA)^{1-\theta} \Big(
            g(\UA)+g(\OA) \Big) \,\nu(d\BA).
\end{align*}
Recall that $Z\neq 0$ by choice of mollifiers. Moreover, at least
one of the constants $A_1$ and $A_2$ is different from zero.
Therefore $B:=(A_1^2+\pi^2 A_2^2)Z$ does not vanish. To conclude
the proof of Proposition~\ref{P:EINS}, we apply the argument above
for the particular choice $g(t) := \LA \chi(t)\RA \zeta(t)$ with
nonnegative $\zeta\in\D(\R)$. As shown in Lemma~\ref{L:CHICON}, the
map $t\mapsto\LA \chi(t) \RA$ is in $C^\alpha(\R)$ for all
$\alpha\in[0,\lambda]$.
\qed


\subsection{Proof of Proposition~\ref{P:ZWEI}}
\label{SS:Z}

We use the notation of Subsection~\ref{SS:E}.

\begin{lemma}
\label{L:EMMA} Let $p\in[1,1/(1-\lambda))$ and let $R\in
W^{1,p}_\LOC(\R)$ be some function. For any distribution
$T\in\{\delta, \PV, H, \log|\cdot|, R\}$ define
$$
    T_\EPS(t) := \int_\R \varphi_\EPS(t-s) T(s) \,ds
    \quad\text{for $(s,\EPS) \in\R\times(0,1)$,}
$$
where $\varphi_\EPS$ is a standard mollifier with $\SPT\varphi_\EPS
\subset[-\EPS,\EPS]$. Then there exists, for any $L>0$, a constant
$C>0$ such that the following estimate holds:
\begin{equation}
    \sup_{\EPS\in(0,1)} \int_0^L t^{\lambda p} |T_\EPS(t)|^p \,dt
        \LS C  \Big( 1+\|R\|^p_{L^\infty(B)} \Big),
\label{E:NEA}
\end{equation}
where $B:=B_{L+2}(0)$. Moreover, as $\EPS\rightarrow 0$ we have
strong convergence
$$
    t_+^\lambda T_\EPS(t) \longrightarrow t_+^\lambda T(t)
    \quad\text{in $L^p_\LOC(\R)$.}
$$
\end{lemma}

{\em Proof.}
Note that $T_\EPS$ is smooth as a function of $\EPS\in(0,1)$. To
establish \eqref{E:NEA} it is therefore sufficient to
consider the behavior as $\EPS\rightarrow 0$. Again we use the
decomposition \eqref{E:COSINE} of $\CI$ into a logarithm and a
smooth function.
\medskip

{\bf Step~1:} We first consider the case of a Dirac measure. We can
estimate
$$
    \bigg| \int_\R \varphi_\EPS(t-s) \delta(s) \,ds \bigg|
        = \varphi_\EPS(t)
        \LS C\EPS^{-1} \mathbf{1}_{[-\EPS,\EPS]}(t),
$$
with $C>0$ some constant depending on $\|\varphi\|_{L^\infty(\R)}$.
Therefore we obtain
\begin{align}
    \int_0^L t^{\lambda p} \bigg| \int_\R \varphi_\EPS(t-s)
            \delta(s) \,ds \bigg|^p \,ds
        & \LS C \EPS^{-p} \int_0^\EPS t^{\lambda p} \,dt
\nonumber\\
        & = C \EPS^{(\lambda-1)p+1} \int_0^1 s^{\lambda p} \,ds,
\label{E:DELMORE}
\end{align}
after a substitution of variables $t=\EPS s$. Since by
assumption $p<1/(1-\lambda)$, the right-hand side of
\eqref{E:DELMORE} converges to zero as $\EPS\rightarrow 0$. This
implies
$$
    t_+^\lambda \bigg( \int_\R \varphi_\EPS(t-s) \delta(s) \,ds
\bigg)
        \longrightarrow 0
    \quad\text{in $L^p_\LOC(\R)$.}
$$
\smallskip

{\bf Step~2:} Now we consider the principal value. Let
$t\in(0,\EPS)$. We decompose
\begin{align}
    & \int_\R \varphi_\EPS(t-s) \PV(s) \,ds
\nonumber\\
    & \quad
        = \int_{-(\EPS-t)}^{\EPS-t} \varphi_\EPS(t-s) \PV(s) \,ds
        + \int_{-(\EPS-t)}^{\EPS+t} \varphi_\EPS(t-s) \PV(s) \,ds.
\label{E:PVMORE}
\end{align}
For the first term we can argue as follows: By symmetry, we have
$$
    \int_{-(\EPS-t)}^{\EPS-t} \varphi_\EPS(t-s) \PV(s) \,ds
        = \int_{-(\EPS-t)}^{\EPS-t} \Big[ \varphi_\EPS(t-s)
            -\varphi_\EPS(t) \Big] \PV(s) \,ds.
$$
Now fix some $\alpha\in(0,1)$. Then we can estimate
\begin{align*}
    \bigg| \int_{-(\EPS-t)}^{\EPS-t} \Big[ \varphi_\EPS(t-s)
            -\varphi_\EPS(t) \Big] \PV(s) \,ds \bigg|
        & \LS \|\varphi_\EPS\|_{C^\alpha(\R)}
            \int_{-(\EPS-t)}^{\EPS-t} |s|^{\alpha-1} \, ds
\\
    & = C \EPS^{-(1+\alpha)} |\EPS-t|^\alpha,
\end{align*}
with $C>0$ some constant depending on $\|\varphi\|_{C^\alpha(\R)}$.
This implies
\begin{align*}
    \int_0^\EPS t^{\lambda p} \bigg| \int_{-(\EPS-t)}^{\EPS-t}
            \varphi_\EPS(t-s) \PV(s) \,ds \bigg|^{p} \,dt
        & \LS C^p \EPS^{-(1+\alpha)p} \int_0^\EPS t^{\lambda p}
            |\EPS-t|^{\alpha p} \,dt
\\
        & = C^p \EPS^{(\lambda-1)p+1} \int_0^1
            s^{\lambda p} |1-s|^{\alpha p} \,ds.
\end{align*}
The right-hand side vanishes as $\EPS\rightarrow 0$. For the second
term in \eqref{E:PVMORE} we find
\begin{align*}
    \left| \int_{\EPS-t}^{\EPS+t} \varphi_\EPS(t-s)
            \PV(s) \,ds \right|
        & \LS \|\varphi_\EPS\|_{L^\infty(\R)}
            \int_{\EPS-t}^{\EPS+t} \frac{ds}{s}
\\
    & = C \EPS^{-1} \bigg| \log\bigg( \frac{\EPS+t}{\EPS-t}
        \bigg) \bigg|,
\end{align*}
with $C>0$ some new constant depending on $\|\varphi\|_{L^\infty
(\R)}$. Therefore
\begin{align*}
    \int_0^\EPS t^{\lambda p} \bigg| \int_{\EPS-t}^{\EPS+t}
            \varphi_\EPS(t-s) \PV(s) \,ds \bigg|^{p} \,dt
        & \LS C^p \EPS^{-p} \int_0^\EPS t^{\lambda p} \bigg|
            \log\bigg( \frac{\EPS+t}{\EPS-t} \bigg) \bigg|^p \,dt
\\
        & = C^p \EPS^{(\lambda-1)p+1} \int_0^1
            s^{\lambda p} \bigg| \log\bigg( \frac{1+s}{1-s} \bigg)
                \bigg|^p \,ds.
\end{align*}
Again the right-hand side converges to zero as $\EPS\rightarrow 0$.
Let now $t\in(\EPS, L)$. Then
$$
    \left| \int_{t-\EPS}^{t+\EPS} \varphi_\EPS(t-s)
            \PV(s) \,ds \right|
        \LS C \EPS^{-1} \bigg| \log\bigg( \frac{t+\EPS}{t-\EPS}
            \bigg) \bigg|,
$$
with $C>0$ some new constant depending on $\|\varphi\|_{L^\infty
(\R)}$. We have
$$
    \sup_{\EPS<t} \EPS^{-1} \bigg| \log\bigg(
            \frac{t+\EPS}{t-\EPS} \bigg) \bigg|
        = \lim_{\EPS\rightarrow 0} \EPS^{-1} \bigg| \log\bigg(
            \frac{t+\EPS}{t-\EPS} \bigg) \bigg|
        = 2t^{-1}.
$$
Therefore we obtain the estimate
\begin{align}
    \int_\EPS^L t^{\lambda p} \bigg| \int_{t-\EPS}^{t+\EPS}
            \varphi_\EPS(t-s) \PV(s) \,ds \bigg|^{p} \,dt
        & \LS (2C)^p \int_\EPS^L t^{(\lambda-1)p} \,dt
\nonumber\\
        & \LS \frac{(2C)^p}{(\lambda-1)p+1} L^{(\lambda-1)p+1}.
\label{E:SUBTLE}
\end{align}
The left-hand side is bounded uniformly in $\EPS$. We conclude that
$$
    t_+^\lambda \bigg( \int_\R \varphi_\EPS(t-s) \PV(s) \,ds \bigg)
        \longrightarrow t_+^{\lambda-1}
    \quad\text{in $L^p_\LOC(\R)$.}
$$
\smallskip

{\bf Step~3:} We now consider the case of a Heaviside function. We
have
$$
    \bigg| \int_\R \varphi_\EPS(t-s) H(s) \,ds \bigg|
        \LS 1.
$$
Therefore we obtain the straightforward estimate
$$
    \int_0^L t^{\lambda p} \bigg| \int_\R \varphi_\EPS(t-s)
            H(s) \,ds \bigg|^p \,ds
        \LS \int_0^L t^{\lambda p} \,dt.
$$
The right-hand side is bounded uniformly in $\EPS$. Moreover, we
have
$$
    t_+^\lambda \bigg( \int_\R \varphi_\EPS(t-s) H(s) \,ds \bigg)
        \longrightarrow t_+^\lambda
    \quad\text{in $L^p_\LOC(\R)$.}
$$
\smallskip

{\bf Step~4:} For the case of a logarithm, we first consider
$t\in(0,\EPS)$. We decompose
\begin{align}
    & \int_\R \varphi_\EPS(t-s) \log|s| \,ds
\nonumber\\
    & \quad
        = \int_{-(\EPS-t)}^0 \varphi_\EPS(t-s) \log|s| \,ds
        + \int_0^{\EPS+t} \varphi_\EPS(t-s) \log|s| \,ds.
\label{E:LOGMORE}
\end{align}
For the first term we can now estimate
\begin{align*}
    \bigg| \int_{-(\EPS-t)}^0 \varphi_\EPS(t-s)
            \log|s| \,ds \bigg|
        & \LS \|\varphi_\EPS\|_{L^\infty(\R)}
            \int_{-(\EPS-t)}^0 |\log|s|| \, ds
\\
    & = C \EPS^{-1} |\EPS-t| \big(1+|\log|\EPS-t||\big)
\end{align*}
with $C>0$ some constant depending on $\|\varphi\|_{L^\infty(\R)}$.
This implies
\begin{align*}
    & \int_0^\EPS t^{\lambda p} \bigg| \int_{-(\EPS-t)}^0
            \varphi_\EPS(t-s) \log|s| \,ds \bigg|^{p} \,dt
\\
    & \quad
        \LS C^p \EPS^{-p} \int_0^\EPS t^{\lambda p}
            |\EPS-t|^p \big(1+|\log|\EPS-t||\big)^p \,dt
\\
    & \quad = C^p \EPS^{\lambda p+1} \big(1+|\log\EPS|\big)^p
        \int_0^1 s^{\lambda p} |1-s|^p \big(1+|\log|1-s||
            \big)^p \,ds.
\end{align*}
The right-hand side vanishes as $\EPS\rightarrow 0$. For the second
term in \eqref{E:LOGMORE} we find
$$
    \bigg| \int_0^{\EPS+t} \varphi_\EPS(t-s)
            \log|s| \,ds \bigg|
        \LS C \EPS^{-1} |\EPS+t| \big(1+|\log|\EPS+t||\big),
$$
which implies the estimate
\begin{align*}
    & \int_0^\EPS t^{\lambda p} \bigg| \int_0^{\EPS+t}
            \varphi_\EPS(t-s) \log|s| \,ds \bigg|^{p} \,dt
\\
    & \quad \LS C^p \EPS^{\lambda p+1} \big(1+|\log\EPS|\big)^p
        \int_0^1 s^{\lambda p} |1+s|^p \big(1+|\log|1+s||
            \big)^p \,ds.
\end{align*}
Again the right-hand side vanishes for $\EPS\rightarrow 0$.
Consider
now $t\in(\EPS,L)$.  Then
\begin{align*}
    & \left| \int_{t-\EPS}^{t+\EPS} \varphi_\EPS(t-s)
            \log|s| \,ds \right|
\\
    & \quad
        \LS C \EPS^{-1} \Big| |t+\EPS|\big( 1+|\log|t+\EPS||
            \big) - |t-\EPS|\big( 1+|\log|t-\EPS|| \big) \Big|,
\end{align*}
with $C>0$ some new constant depending on $\|\varphi\|_{L^\infty
(\R)}$. We have
\begin{align*}
    & \sup_{\EPS<t} \EPS^{-1} \Big| |t+\EPS|\big( 1+|\log|t+\EPS||
            \big) - |t-\EPS|\big( 1+|\log|t-\EPS|| \big) \Big|
\\
    & \quad
        = \lim_{\EPS\rightarrow 0} \EPS^{-1} \Big| |t+\EPS|
            \big( 1+|\log|t+\EPS|| \big) - |t-\EPS|\big(
                1+|\log|t-\EPS|| \big) \Big|
        = 2|\log|t||.
\end{align*}
Therefore we obtain the estimate
$$
    \int_\EPS^L t^{\lambda p} \bigg| \int_{t-\EPS}^{t+\EPS}
            \varphi_\EPS(t-s) \log|s| \,ds \bigg|^{p} \,dt
        \LS (2C)^p \int_\EPS^L t^{\lambda p} |\log|t||^p \,dt.
$$
The right-hand side is bounded uniformly in $\EPS$. We obtain
$$
    t_+^\lambda \bigg( \int_\R \varphi_\EPS(t-s) \log|s| \,ds
        \bigg) \longrightarrow t_+^\lambda \log|t|
    \quad\text{in $L^p_\LOC(\R)$.}
$$
\smallskip

{\bf Step~5:} Finally, let us consider the case of a function $R\in
W^{1,p}_\LOC(\R)$. By Sobolev embedding theorems, the function
$R\in
C^\alpha(\R)$ for some $\alpha\in[0,\lambda)$. We have
$$
    \int_0^L t^{\lambda p} \bigg| \int_\R \varphi_\EPS(t-s)
            R(s) \,ds \bigg|^p \,ds
        \LS \|R\|^p_{L^\infty(B)}
            \int_0^L t^{\lambda p} \,dt,
$$
using Minkowski inequality. The convergence
$$
    t_+^\lambda \bigg( \int_\R \varphi_\EPS(t-s) R(s) \,ds
        \bigg) \longrightarrow t_+^\lambda R(t)
    \quad\text{in $L^p_\LOC(\R)$}
$$
follows from well-known results on mollification of
$L^p_\LOC$-functions.
\qed

\begin{remark}\label{R:EMARK}
A careful inspection of the previous proof shows that the statement
of Lemma~\ref{L:EMMA} is still true for $T\in\{H, \CI, R\}$ and
$t_+^{\lambda-1}$. We have
$$
    \sup_{\EPS\in(0,1)}
            \int_0^L t^{(\lambda-1)p} |T_\EPS(t)|^p \,dt
        \LS C \Big( 1+\|R\|^p_{L^\infty(B)} \Big)
$$
for some constant $C>0$ depending on $L$, and the strong
convergence
$$
    t_+^{\lambda-1} T_\EPS(t) \longrightarrow t_+^{\lambda-1} T(t)
    \quad\text{in $L^p_\LOC(\R)$.}
$$
For $T\in\{\delta,\PV\}$ and $t_+^{\lambda-1}$ we obtain the bound
$$
    \sup_{\EPS\in(0,1)} \EPS^p
            \int_0^L t^{(\lambda-1)p} |T_\EPS(t)|^p \,dt
        \LS C
$$
for some $C>0$. Note the extra factor $\EPS^p$ needed here to
control the integral. Again the necessary estimates can be adapted
easily. We have
\begin{align}
    \EPS^p \int_\EPS^L t^{\lambda p} \bigg| \int_{t-\EPS}^{t+\EPS}
            \varphi_\EPS(t-s) \PV(s) \,ds \bigg|^{p} \,dt
        & \LS \EPS^p (2C)^p \int_\EPS^L t^{(\lambda-2)p} \,dt
\nonumber\\
        & \LS \frac{(2C)^p}{|(\lambda-2)p+1|} \EPS^{(\lambda-1)p+1}
\label{E:SUBTLE2}
\end{align}
instead of \eqref{E:SUBTLE}. The right-hand side of
\eqref{E:SUBTLE2} converges to zero as $\EPS\rightarrow 0$.
\end{remark}

\begin{lemma}
\label{L:ALMOST}
Let $f(s)=(1-s^2)^\lambda_+$ for all $s\in\R$. Fix
some $p\in[1, 1/(1-\lambda))$ and a standard mollifier
$\varphi_\EPS$ such that $\SPT \varphi_\EPS\subset[-\EPS,\EPS]$.
Then we have
\begin{equation}
\begin{aligned}
    \sup_{\EPS\in(0,1)} \bigg\|f(t)\bigg(\int_\R
        \varphi_\EPS(t-s)\,\DS
            f(s)\,ds \bigg) \bigg\|_{W^{1,p}(\R)}
        &\LS C \Big( 1+\|r\|_{L^\infty(\R)} \Big),
\\
    \sup_{\EPS\in(0,1)} \bigg\|f(t)\bigg(\int_\R (t-s)
        \varphi_\EPS(t-s)
            \,\DC f(s)\,ds \bigg) \bigg\|_{W^{1,p}(\R)}
        &\LS C \Big( 1+\|q\|_{L^\infty(\R)} \Big),
\end{aligned}
\label{E:UNEA}
\end{equation}
with $C>0$ some constant. Moreover, we find
\begin{equation}
\left.\begin{aligned}
    f(t) \bigg( \int_\R \varphi_\EPS(t-s) \,\DS f(s) \,ds \bigg)
        &\longrightarrow f(t)\, \DS f(t)
\\
    f(t) \bigg( \int_\R (t-s)\varphi_\EPS(t-s) \,\DC f(s) \,ds
        \bigg) &\longrightarrow 0
\end{aligned}\right\}
    \quad\text{in $W^{1,p}(\R)$}
\label{E:STATE}
\end{equation}
as $\EPS\rightarrow 0$. This implies strong convergence in
$C^\alpha(\R)$, for some $\alpha\in[0,\lambda)$.
\end{lemma}

{\em Proof.} Note first that by Proposition~\ref{P:EXP}, the
derivative $\DS f$ contains Heavi\-side functions, logarithms and a
remainder in $W^{1,p}_\LOC(\R)$. We have
\begin{align*}
    & \frac{d}{dt} \bigg\{ f(t) \bigg( \int_\R \varphi_\EPS(t-s)
            \,\DS f(s) \,ds \bigg) \bigg\}
\\
    & \quad
        = \frac{df(t)}{dt} \bigg( \int_\R \varphi_\EPS(t-s)
            \,\DS f(s) \,ds \bigg)
        + f(t) \bigg( \int_\R \varphi_\EPS(t-s) \,\DC f(s)
            \,ds \bigg)
\end{align*}
for a.e.\ $t\in\R$, where we used \eqref{E:RFV1}. The derivative of
$f(t)$ blows up like $|1-|t||_+^{\lambda-1}$ as $|t|\rightarrow 1$.
We apply Lemma~\ref{L:EMMA} and Remark~\ref{R:EMARK} to obtain
$$
    \frac{d}{dt} \bigg\{ f(t) \bigg( \int_\R \varphi_\EPS(t-s)
            \,\DS f(s) \,ds \bigg) \bigg\}
        \longrightarrow \frac{df(t)}{dt}\, \DS f(t) + f(t)\, \DC
f(t)
    \quad\text{in $L^p(\R)$}
$$
as $\EPS\rightarrow 0$. The first statement in \eqref{E:STATE}
follows. Similarly, we write
\begin{align}
    & \frac{d}{dt} \bigg\{ f(t) \bigg( \int_\R
(t-s)\varphi_\EPS(t-s)
            \,\DC f(s) \,ds \bigg) \bigg\}
\nonumber\\
    & \quad
        = \EPS \frac{df(t)}{dt} \bigg( \int_\R \psi_\EPS(t-s)
            \,\DC f(s) \,ds \bigg)
        + f(t) \bigg( \int_\R (\partial_t\psi)_\EPS(t-s) \,\DC f(s)
            \,ds \bigg),
\label{E:MORE2}
\end{align}
with $\psi(t) := t\varphi(t)$ and $\psi_\EPS(t) := \EPS^{-1}
\psi(t/\EPS)$ for all $(s,\EPS)\in\R\times(0,1)$. We apply
Lemma~\ref{L:EMMA} and Remark~\ref{R:EMARK} to obtain the second
bound in \eqref{E:UNEA} and convergence in $L^p(\R)$ as
$\EPS\rightarrow 0$. Note that the extra factor $\EPS$ causes the
first term on the right-hand side of \eqref{E:MORE2} to vanish. For
the second term we apply the dominated convergence theorem: Since
$\partial_t \psi$ has zero mean, we have pointwise convergence to
zero almost everywhere. We conclude that
$$
    \frac{d}{dt} \bigg\{ f(t) \bigg( \int_\R (t-s)\varphi_\EPS(t-s)
            \,\DC f(s) \,ds \bigg) \bigg\}
        \longrightarrow 0
    \quad\text{in $L^p(\R)$}
$$
as $\EPS\rightarrow 0$, which implies the second statement in
\eqref{E:STATE}.
\qed

{\em Proof of Proposition~\ref{P:ZWEI}.} Using \eqref{E:CHISIG} and
\eqref{E:RFV2} we find the identity
\begin{align}
    & \chi(t|\BA)\, \DC\sigma(s|\BA)
        -\sigma(t|\BA)\, \DC\chi(s|\BA)
\nonumber\\
    & \vphantom{\Big[}\quad
        = \theta (t-s)\, \chi(t|\BA)\, \DC\chi(s|\BA)
        + \theta(\lambda+1) \chi(t|\BA)\, \DS\chi(s|\BA),
\label{E:CSMSC3}
\end{align}
which holds distributionally in $(s,s')\in\R\times\R$ for all
$\BA\in\CC$. Let us consider the first term on the right-hand side.
We fix some $\BA\in\CC$ and integrate against the mollifier
$\varphi_\EPS(t-s)$. We apply Lemmas~\ref{L:EMMA} and \ref{L:ALMOST}
and obtain that
$$
    \bigg\|\chi(t|\BA)\bigg(\int_\R \varphi_\EPS(t-s)
            \,\DS\chi(s|\BA)\,ds \bigg) \bigg\|_{W^{1,p}(K)}
        \LS C \rho(\BA)^{3\theta\lambda}
$$
for all $K\subset\R$ compact, with $C>0$ depending on
$K$ and $\|r\|_{L^\infty(\R)}$. Recall that $0<3\theta\lambda
<\gamma+1$ for $\gamma\in(1,3)$. We can integrate against $\nu$ to
get
$$
    \bigg\|\bigg\LA \chi(t)\bigg(\int_\R \varphi_\EPS(t-s)
            \,\DS\chi(s)\,ds \bigg) \bigg\RA\bigg\|_{W^{1,p}(K)}
        \LS C \int_\CC W(\BA) \,\nu(d\BA),
$$
which is finite by assumption on $\nu$. Sending $\EPS\rightarrow
0$,
we obtain
\begin{equation}
    \bigg\LA \chi(t)\bigg(\int_\R \varphi_\EPS(t-s)
            \,\DS\chi(s)\,ds \bigg) \bigg\RA
        \longrightarrow \big\LA \chi(t) \,\DS\chi(t) \big\RA
    \quad\text{locally in $C^\alpha(\R)$,}
\label{E:EEN}
\end{equation}
for some $\alpha\in(0,\lambda)$. We used Lemma~\ref{L:ALMOST} and
Sobolev embedding. Similarly
\begin{align*}
    & \bigg\|\chi(t|\BA)\bigg(\int_\R (t-s)\varphi_\EPS(t-s)
            \,\DC\chi(s|\BA)\,ds \bigg) \bigg\|_{W^{1,p}(K)}
\\
    & \quad
        \LS C \rho(\BA)^{3\theta\lambda}
                \Big( 1+\rho(\BA)^{-\theta} \Big)
                \Big( 1 + |\log\rho(\BA)| \Big),
\end{align*}
with $C>0$ some constant. Since $0<(3\lambda-1)\theta<\gamma+1$
for $\gamma\in(1,3)$, we get
$$
    \bigg\|\bigg\LA \chi(t)\bigg(\int_\R (t-s)
            \varphi_\EPS(t-s) \,\DS\chi(s)\,ds \bigg)
                \bigg\RA\bigg\|_{W^{1,p}(K)}
        \LS C \int_\CC W(\BA) \,\nu(d\BA).
$$
Sending $\EPS\rightarrow 0$, we obtain that
\begin{equation}
    \bigg\LA \chi(t)\bigg(\int_\R (t-s)\varphi_\EPS(t-s)
            \,\DC\chi(s)\,ds \bigg) \bigg\RA
        \longrightarrow 0
    \quad\text{locally in $C^\alpha(\R)$,}
\label{E:TWE}
\end{equation}
as follows from Lemma~\ref{L:ALMOST} and Sobolev embedding.
Therefore
\begin{equation}
    \Big\LA \chi(t) \DC\sigma_\EPS(t)-\sigma(t)\, \DC\chi_\EPS(t)
        \Big\RA \longrightarrow \theta(\lambda+1) \big\LA \chi(t)
            \,\DS\chi(t) \big\RA
    \quad\text{locally in $C^\alpha(\R)$.}
\label{E:FINA}
\end{equation}
Note that \eqref{E:EEN} and \eqref{E:TWE} are independent of the
choice of mollifier: we can use $\varphi'_\EPS(t-s)$
instead (see the beginning of Subsection~\ref{SS:RED} for the
definition) and obtain the analogous convergence as in
\eqref{E:FINA}, with the same limit.

To conclude the proof of Proposition~\ref{P:ZWEI} it is now
sufficient to notice that
\begin{equation}
    \big\LA \DC\chi'_\EPS(t) \big\RA
        \WEAK \big\LA \DC\chi(t) \big\RA
    \quad\text{weakly-$\star$ in $\big( C^\alpha_c(\R)
\big)^{\!*}$}
\label{E:CONDUA}
\end{equation}
(the dual of the space of H\"{o}lder continuous functions with
compact support). Recall that the fractional derivative
$\DC\chi(\cdot|\BA)$ contains only Dirac measures, principal value
operators, and locally integrable functions (see \eqref{E:DLP}). It
stays bounded uniformly as $\rho(\BA) \rightarrow 0$ since $\lambda
\GS 1$ if $\gamma\in(1,5/3]$, and grows at most linearly for
$\rho(\BA)$ large. Recall that if $\gamma=5/3$, then the constant
$A_4$ in \eqref{E:DLP} vanishes, so the logarithmic term does not
matter. We can now integrate $\DC\chi(\cdot|\BA)$ against $\nu$,
and
then \eqref{E:CONDUA} follows. The same convergence holds if we use
the mollifier $\varphi_\EPS(t-s)$ instead.

For any test function $\zeta\in\D(\R)$ we therefore obtain
\begin{align*}
    & \lim_{\EPS\rightarrow 0} \int_\R \Big\LA
                 \chi(t) \DC\sigma'_\EPS(t)
                -\sigma(t) \DC\chi'_\EPS(t) \Big\RA
            \big\LA\DC\chi_\EPS(t)\big\RA \zeta(t) \,dt
\\
    & \quad
        = \lim_{\EPS\rightarrow 0} \int_\R \Big\LA
                 \chi(t) \DC\sigma_\EPS(t)
                -\sigma(t) \DC\chi_\EPS(t) \Big\RA
            \big\LA\DC\chi'_\EPS(t)\big\RA \zeta(t) \,dt
\\
    & \quad = \theta(\lambda+1) \int_\R \big\LA \chi(t)
            \DS\chi(t) \big\RA \big\LA\DC\chi(t) \big\RA
                \zeta(t) \,dt.
\end{align*}
This completes the proof of the proposition.
\qed


\appendix
\section{Propagation of equi-integrability}

For nozzle flows with $A$ constant, the proof of
Proposition~\ref{P:NONL} can also be based on the following lemma,
which shows that for entropy solutions of the isentropic Euler
equations, equi-integrability of the total energy is
``propagated.''
We complement assumptions (i)--(iv) of Section~\ref{SS:APPROX} by
requiring that
\begin{enumerate}
\item[(v)] the sequence $(\OVR^n,\OVU^n)$ vanishes uniformly in the
large
in the sense that for each $\EPS>0$ there exists a compact subset
$K\subset\R$ with
$$
    \sup_n \int_{\R\setminus K}
        \Big( \HA\OVR^n (\OVU^n)^2
            +U(\OVR^n) \Big) A^n \, dx
    \LS \EPS.
$$
\end{enumerate}
Under this assumption, \eqref{E:CLA} of Lemma~\ref{L:CONID} can be
improved: With the notation used there, we have that for all
$\EPS>0$ there exist $N,R>0$ such that
\begin{equation}
    \sup_{n\GS N} \iint_{\R\times\R} s^2\Phi_R(s)\,
        \chi(s|\OVZ^n) \,ds \,dx \LS \EPS.
\label{E:CLA2}
\end{equation}

Then we have the following result.

\begin{lemma}\label{L:EQUIINT} Choose a test function
$\varphi\in\D(\R)$ with $0\LS\varphi\LS 1$, such that
$\varphi(s)=1$ for $|s|\LS 1$ and $\varphi(s)=0$ for $|s|\GS 2$.
Define $\varphi_R := \varphi(\cdot /R)$ and $\Phi_R :=
1-\varphi_R$. For all $T>0$ and all $\EPS>0$ there exist $R,N>0$
such that
\begin{gather}
    \sup_{n\GS N}
        \iint_{[0,T]\times\R} \int_\R s^2 \Phi_R(s)\,
            \chi(s|\BZ^n) \, ds \, dx\,dt
        \LS \EPS,
\label{E:RESU1}
\\
    \sup_{n\GS N}
        \iint_{[0,T]\times\R} \int_\R |s| \Phi_R(s)\,
            |\sigma(s|\BZ^n)| \,ds \, dx\,dt
        \LS \EPS.
\label{E:RESU2}
\end{gather}
\end{lemma}

{\em Proof.}
By \eqref{E:CLA2}, there exist $R,N>0$ such that
\begin{equation}
    \sup_{n\GS N} \iint_{\R\times\R} 2s^2\Phi_{R/2}(s)\,
        \chi(s|\OVZ^n) \,ds \,dx \LS \EPS/T.
\label{E:INTMO}
\end{equation}
For this $R$ let $\psi(s) :=  2(s^2-R^2)\, \IND_{\{|s|\GS R\}}$ for
all $s\in\R$. Since $\psi$ is convex we can use this weight function
in the entropy inequality \eqref{E:ENTROPY} and obtain
\begin{equation}
    \ESUP_{t\GS 0} \iint_{\R\times\R}
            \psi(s)\, \chi(s|\BZ^n(t,x)) \,ds \,dx
        \LS \int_{\R\times\R} \psi(s)\, \chi(s|\OVZ^n)
                \,ds \,dx
\label{E:ENT}
\end{equation}
for all $n$. On the other hand, we have the following estimate:
$$
    s^2\Phi_R(s)
        \LS \psi(s)
        \LS 2s^2\Phi_{R/2}(s)
    \quad\text{for all $s\in\R$.}
$$
Combining this with \eqref{E:INTMO} and \eqref{E:ENT}, we find that
for all $n\GS N$
$$
    \ESUP_{t\GS 0} \iint_{\R\times\R} s^2\Phi_R(s)\,
        \chi(s|\BZ^n(t,x)) \,ds \,dx \LS \EPS/T,
$$
and integrating over $[0,T]$ we obtain \eqref{E:RESU1}.

To derive \eqref{E:RESU2}, we use the estimate
\begin{align*}
    & \iint_{\R^2} |s|\Phi_R(s)\,
        |\sigma(s|\BZ^n(t,x))| \,ds \,dx
    \LS \theta \iint_{\R^2} s^2\Phi_R(s)\,
        \chi(s|\BZ^n(t,x)) \,ds \,dx
\\
    & \qquad
        +(1-\theta) \bigg( \int_\R \big( \rho^n(u^n)^2 \big)(t,x)
            \,dx \bigg)^{\!1/2}
\\
    & \hspace{10em}
        \bigg( \iint_{\R\times\R} s^2\Phi_R(s)\,
            \chi(s|\BZ^n(t,x)) \,ds \,dx \bigg)^{\!1/2}
\end{align*}
for almost every $t$. The kinetic energy is uniformly bounded by
\eqref{E:ENTRINEQ}.
\qed


\section*{Acknowledgments} 
The first author (P.G.L.) was supported by the A.N.R. Grant
06-2-134423: {\em Mathematical methods in general relativity}
(MATH-GR) and the Centre National de la Recherche Scientifique
(CNRS). The second author (M.W.) acknowledges partial support by the
European network grant  HRPN-CT-2002-00282: {\em Hyperbolic and
kinetic equations,} and by the research project
``Sonderforschungsbereich 611'' {\em Singular phenomena and scaling in
mathematical models} at Bonn University.


\end{document}